\documentclass[a4paper]{article}
\pdfoutput=1
\usepackage[utf8]{inputenc}
\usepackage{amsmath, amssymb, amsthm}
\usepackage{graphicx}
\usepackage[colorlinks=true,allcolors=black]{hyperref}

\usepackage{caption}
\usepackage{subcaption}

\usepackage{dsfont}
\usepackage{overpic}

\usepackage{pgfplots}
\pgfplotsset{compat=1.13}

\usepackage{tikz,tikz-3dplot}
\usetikzlibrary{external} 
\renewcommand{\d}{{\operatorname{d}}}

\newcount\Comments  
\newcommand{\kibitz}[2]{\ifnum\Comments=1\textcolor{#1}{#2}\fi}

\title{Simplex Stochastic Collocation for \\Piecewise Smooth Functions with Kinks}
\usepackage{authblk}
\author[1]{Barbara Fuchs}
\author[1,2,3]{Jochen Garcke}
\affil[1]{Fraunhofer SCAI, Sankt Augustin, Germany}
\affil[2]{Fraunhofer Center for Machine Learning, Sankt Augustin, Germany }
\affil[3]{Institut für Numerische Simulation, Universität Bonn, Germany}
\begin{document}
\maketitle

\paragraph{Abstract}
Most approximation methods in high dimensions exploit smoothness of the function being approximated. 
These methods provide poor convergence results for non-smooth functions with kinks. For example, 
such kinks can arise in the uncertainty quantification of  quantities of interest for gas networks. 
This is due to the regulation of the gas flow, pressure, or temperature. But, one can exploit that 
for each sample in the parameter space it is known if a regulator was active or not, which can be 
obtained from the result of the corresponding numerical solution. 
This information can be exploited in a stochastic collocation method. We approximate the function 
separately on each smooth region by polynomial interpolation and obtain an approximation to the 
kink. 
Note that we do not need information about the exact location of kinks, but only an indicator 
assigning each sample point to its smooth region. We obtain a global order of convergence of 
$(p+1)/d$, where $p$ is the degree of the employed polynomials and $d$ the dimension of the 
parameter space.

\section{Introduction}
In many applications in engineering and science the input data of meta models or simulations is 
uncertain. These uncertainties can arise, e.g.\ in the geometry, boundary conditions, or model 
coefficients. One is often interested in how these uncertainties influence some specific output 
variables also called quantities of interest (QoI). For such an uncertainty quantification (UQ) we 
need methods to approximate and integrate high dimensional functions. In the case of smooth 
functions there are several methods such as (adaptive) sparse grids, Galerkin methods, Polynomial 
chaos expansion or quasi Monte Carlo methods~\cite{Smith2014}. 
For discontinuous functions methods exist such as adaptive sparse grids \cite{jakeman2011}, Voronoi 
piecewise surrogate models (VPS) \cite{rushdi2017} or simplex stochastic collocation (SSC) 
\cite{witteveen2012_1, witteveen2012_2, witteveen2013}. The ideas behind VPS and SSC are similar, in 
both cases the function is locally approximated by piecewise polynomials either on Voronoi cells or 
on simplices resulting from a Delaunay triangulation. In VPS a jump in the function is detected if 
the difference in the function values between neighboring cells exceeds a user defined threshold, 
while SSC detects a jump not directly but by observing the resulting oscillations. Other alternative 
approaches for handling discontinuities include enriching the polynomial approximation basis, which 
generally requires some a priori knowledge of the discontinuity, domain decomposition, also known as 
multi-element approximations in this context, or discontinuity detection algorithms, see 
e.g.~\cite{Colombo2018,Sargsyan2012} for current references. Note that non-smooth functions with 
kinks can be smoothed by integration \cite{griebel2010, griebel2015} over one dimension if the 
location of the kink is known.

In the simulation of gas networks the solution functions are continuous, but not globally 
differentiable due to human intervention through the use of control valves, compressors, or 
heaters. Kinks in a function arise at hyper-surfaces where the function is not continuously 
differentiable. The idea of simplex stochastic collocation~\cite{witteveen2012_1, witteveen2012_2, 
witteveen2013} is to approximate a function $f$ by a piecewise polynomial interpolation on 
simplices. Since polynomial interpolation gets oscillatory near discontinuities, one ensures that 
the approximation is local extremum conserving, i.e.\ maximum and minimum of the approximation in 
any simplex must be attained at its vertices, otherwise the polynomial degree is decreased by 
one~\cite{witteveen2012_1}. This condition results in a fine discretization near discontinuities and 
a coarser discretization at smooth regions. We evaluated the original approach for functions with 
kinks, but were not able to reach the desired convergence rates, i.e.\ by increasing the polynomial 
degree the approximation of a kink could not be improved. 

Based on the original simplex stochastic collocation, we introduce a new approach by taking 
advantage of additional knowledge. In particular, we assume to know on which side of the kink a 
specific collocation point is situated. 
This enables us to approximate the function on each side of the kink separately. In doing so, we can 
improve the convergence rate significantly by not wasting sampling points near the kink. This 
assumption is motivated by the uncertainty quantification for gas networks. Although we have no 
information regarding the location of a kink, we know which elements of the gas network cause kinks. 
After simulating the gas flow for a specific combination of uncertain parameters, we know whether 
the kink inducing elements are active or not. In the case of a control valve we only need to check 
if the outgoing pressure lies below the preset pressure $p_\text{set}$ or equals it. 

The paper is organized as follows. In the second section we introduce the SSC method in general and 
our modifications for piecewise smooth functions with kinks. In addition, we discuss where to sample 
a new point in order to refine a simplex, how to estimate the error, and if it is possible to refine 
multiple simplices at once. In the third section we quantify the uncertainty in a particular node of 
a gas network caused by uncertain input data by applying SSC to calculate the expected pressure.

\section{Simplex Stochastic Collocation}
We now introduce the approach of simplex stochastic collocation following \cite{witteveen2012_1, 
witteveen2012_2, witteveen2013}. Let $\Omega=[0,1]^d$ and $f: \Omega \rightarrow \mathbb{R}$ be a 
continuous function. We first discuss the Delaunay triangulation of a given set of $n$ uniformly 
distributed sampling points $\mathbf{x}_i$, which divides the parameter space $\Omega$ into $m$ 
disjoint simplices $T_j$, before considering refinement strategies. Note that the sampling points 
always include the corners of $\Omega$. Each simplex $T_j$ is defined by its $d+1$ vertices 
$\mathbf{x}_{i_{j,l}}$, with $i_{j,l} \in \{1,\ldots,n\}$ and $l\in \{0, \ldots, d\}$. 

\subsection{The Original SSC}
\begin{figure}
	\centering
	\begin{subfigure}{0.32\textwidth}
		\centering
		\includegraphics{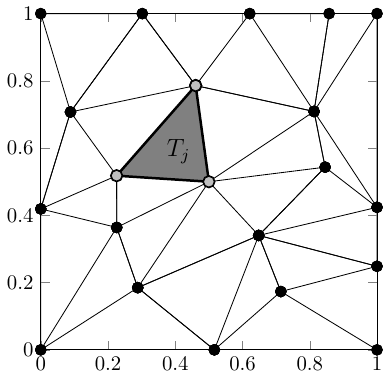}
		\subcaption{$p_j=1$, 3-nn}
	\end{subfigure}
	\hfill
	\begin{subfigure}{0.32\textwidth}
		\centering
		\includegraphics{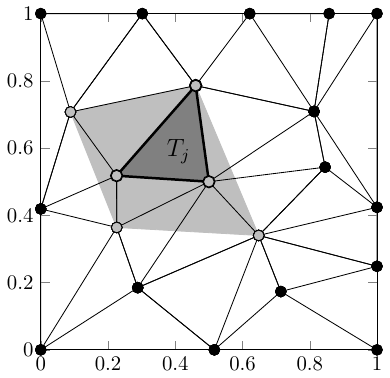}
		\subcaption{$p_j=2$, 6-nn}
	\end{subfigure}
	\hfill
	\begin{subfigure}{0.32\textwidth}
		\centering
		\includegraphics{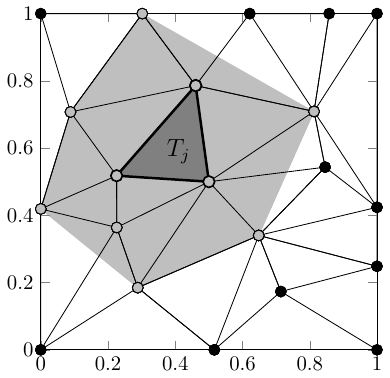}
		\subcaption{$p_j=3$, 10-nn}
	\end{subfigure}	
	\caption{Shown are the Delaunay triangulation of $n=20$ sampling points and the nearest neighbor 
	stencils $S_j$ (light gray) for simplex $T_j$ (dark gray) for polynomial 	degrees $p_j=1,2,3$.}
	\label{stencil}
\end{figure}

Let $f\in\mathcal{C}^0(\Omega)$ be a continuous function that we approximate by $m$ piecewise 
polynomial functions $g_j(\mathbf{x})$ defined on simplex $T_j$
\begin{align*}
	f(\mathbf{x}) \approx \sum_{i=1}^m g_j(\mathbf{x}) \;\mathds{1}_{T_j}.
\end{align*}
The polynomials $g_j$ are defined as
\begin{align*}
	g_j(\mathbf{x}) = \sum_{k=1}^{N_j} c_{j,k}\psi_{j,k}(\mathbf{x}),
\end{align*}
where $\psi_{j,k}$ are some appropriate basis polynomials, $c_{j,k}$ the corresponding coefficients, 
and $N_j=(d+p_j)!/(d!p_j!)$ the number of degrees of freedom, with $p_j \leq p_\text{max}$ the local 
polynomial degree. Note that in our numerical experiments we use the monomial basis. The polynomial 
approximation $g_j(\mathbf{x})$ in $T_j$ is constructed by interpolating $f(\mathbf{x})$ in a 
stencil 
\begin{align*}
	S_j = \{\mathbf{x}_{i_{j,0}}, \ldots, \mathbf{x}_{i_{j,N_j}}\}
\end{align*}
consisting of $N_j$ points out of the sampling points $\mathbf{x}_i$. These points are chosen to be 
the nearest neighbors to simplex $T_j$ based on the Euclidean distance to its center of mass. Since 
in the case of long and flat simplices not necessarily all of its vertices belong to the set of 
nearest neighbors, we always include the $d+1$ simplex vertices in $S_j$. Thus, we ensure that our 
approximation is exact at all sampling points. See Figure \ref{stencil} for different nearest 
neighbor stencils of simplex $T_j$ corresponding to polynomial degrees $p_j=1,2,3$. If the 
interpolation problem is not uniquely solvable we reduce the polynomial degree $p_j$ successively by 
one until the solution is unique. To avoid oscillations in an approximation $g_j(\mathbf{x})$ near a 
discontinuity, the local polynomial degree $p_j$ is also reduced by one if the approximation is not 
local extremum conserving (LEC), i.e.\ if it does not hold that
\begin{align}
	\min_{\mathbf{x} \in T_j} g_j(\mathbf{x}) = \min_{\mathbf{x}_i \in T_j} f(\mathbf{x}_i)
	\quad \wedge \quad
	\max_{\mathbf{x} \in T_j} g_j(\mathbf{x}) = \max_{\mathbf{x}_i \in T_j} f(\mathbf{x}_i).
	\label{LEC}
\end{align}
Note that the polynomial degree will be at least one. This holds because the linear interpolation 
problem on a simplex is always uniquely solvable and the resulting interpolation is always local 
extremum conserving. Since the approximation in one single simplex is independent from all other 
simplices, the resulting global approximation is not even continuous across the simplices' facets, 
except for linear polynomials. 

\subsubsection{The Theoretical Convergence Rate}
For smooth functions $f\in\mathcal{C}^{p+1}$ and uniformly distributed sampling points we can 
locally estimate the approximation error. Let $\{\mathbf{x}_\alpha\}_{|\alpha|\leq p}$ denote the 
interpolation points with multi-index $\alpha=(\alpha_1, \alpha_2, \ldots, \alpha_d) \allowbreak \in 
\mathbb{N}_0^d$. The classic estimation \cite{sauer1995} for the error in the $d$-dimensional point 
$\mathbf{x} = (x^{(1)}, \ldots, x^{(d)})$ between the function $f(\mathbf{x})$ and its Lagrange 
interpolation $L_p f(\mathbf{x})$ of degree $p$ reads
\begin{align}
	\left|L_p f(\mathbf{x}) - f(\mathbf{x})\right| 
	\leq \hspace{-0.5em} \sum_{|\alpha|=p+1} \frac{1}{\alpha!} \left\| 
\frac{\partial^{p+1}f}{\partial 
	\mathbf{x}^\alpha} \right\|_\infty & \prod_{\gamma_1=1}^{\alpha_1} \hspace{-0.2em} \left(x^{(1)} 
- 	x^{(1)}_{(\gamma_1-1,\alpha_2, \ldots, \alpha_d)}\right) \notag \\
	\quad \cdots & \prod_{\gamma_d=1}^{\alpha_d} \hspace{-0.2em} \left(x^{(d)} 
	- x^{(d)}_{(\alpha_1,\alpha_2,\ldots, \gamma_d-1)}\right).
	\label{err:lagrange}
\end{align}
In the $i$-th product the $i$-th entry of $\alpha$ is replaced by $\gamma_i-1$. For $n$ uniformly 
distributed random points in $\Omega$ the expected distance between two of them is of order 
$\mathcal{O}(n^{-1/d})$. Because each summand consists of $p+1$ factors, each summand is of order 
$\mathcal{O}(n^{-(p+1)/d})$. Thereby we can estimate the products in (\ref{err:lagrange}) and obtain
\begin{align}
	\left|L_p f(\mathbf{x}) - f(\mathbf{x})\right| \leq C \cdot n^{-(p+1)/d} \sum_{|\alpha|=p+1} 
	\frac{1}{\alpha!} \left\| \frac{\partial^{p+1}f}{\partial \mathbf{x}^\alpha} \right\|_\infty .
	\label{eq:convergence}
\end{align}
Thus the Lagrange interpolation $L_p f$ converges pointwise with order $(p+1)/d$ against the 
function $f$ if the partial derivatives are bounded. Because the order of convergence depends on 
the dimension we need to increase the polynomial degree with increasing dimension to obtain a 
constant order of convergence. The error estimate (\ref{eq:convergence}) holds true for any simplex 
$T_j$ and corresponding approximation $g_j(\mathbf{x})$. Note that for functions $f\in 
\mathcal{C}^0(\Omega)$ with kinks, i.e.\ functions that are continuous but not continuously 
differentiable, we cannot estimate the error with (\ref{eq:convergence}) or expect an order of 
convergence of $(p+1)/d$, as $f\notin \mathcal{C}^{p+1}(\Omega)$. This motivates the following 
modification of the original approach. 

\subsection{The Improved SSC}
Let $f\in\mathcal{C}^0(\Omega)$ be a function with kinks. We say a function $f:[0,1]^d \to 
\mathbb{R}$ has a kink at the $(d-1)$-dimensional hyper-surface $K\subset\Omega$ if for all 
$\mathbf{x} \in K$ the function $f(\mathbf{x})$ is not continuously differentiable. In $d=2$ 
dimensions, the kink locations are lines and can be arbitrarily shaped, they can be straight, curved 
or closed lines, and they can also intersect. In $d=3$ dimensions, the kink locations are surfaces. 
We applied the approaches from \cite{witteveen2012_1, witteveen2012_2, witteveen2013} to functions 
with kinks, but were not able to reach the desired convergence rates, as can be seen in our 
numerical experiments in sections~\ref{sec:num_test} and \ref{sec:gas_net_results}. Therefore, we 
developed an improved SSC, which we introduce in the following.

\begin{figure}
\begin{center}
	\begin{subfigure}{0.45\textwidth}
		\centering
		\includegraphics{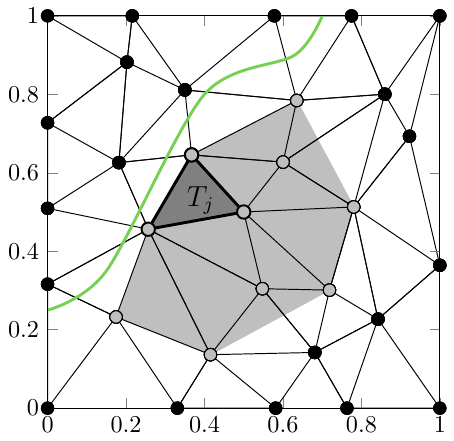}
		\subcaption{Simplex without a kink.}
	\label{improved_stencil_1}
	\end{subfigure}
	\begin{subfigure}{0.45\textwidth}
		\centering
		\includegraphics{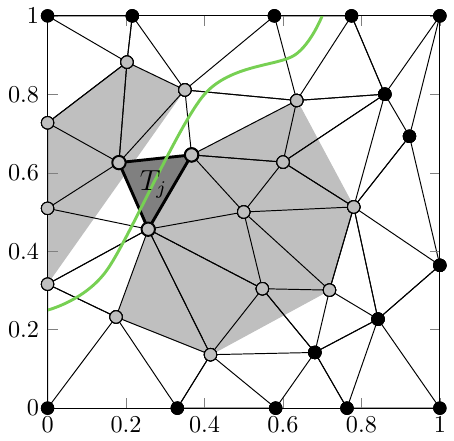}
		\subcaption{Simplex with a kink.}
	\label{improved_stencil_2}
	\end{subfigure}
\end{center}
	\caption{The improved nearest-neighbor stencils (light gray) for simplex $T_j$ (dark gray). The 
	domain $\Omega$ is divided by a kink (green) into two subdomains $\Omega_1$ (left) and $\Omega_2$ 
	(right). (a) shown is a stencil for a simplex without a kink 	inside and completely lying in 
	$\Omega_2$. (b) shown are two stencils for a simplex 	with a kink inside and lying in $\Omega_1$ 
	as well as in $\Omega_2$.}
	\label{improved_stencil}
\end{figure}

Generally, kinks divide the parameter space $\Omega$ into disjoint subdomains $\Omega_k$ with 
$\bigcup_k \Omega_k = \Omega$. Suppose $f\in\mathcal{C}^{p+1}(\Omega_k)$ is smooth for all $k$, and 
that we have for each sampling point $\mathbf{x}_i$ the information to which $\Omega_k$ it belongs. 
The last assumption is motivated by our application of gas networks, where one knows if a regulator 
is active or not, which influences the locations of kinks. There are two different cases for our 
modification:

\paragraph{Case 1.} The function $f(\mathbf{x})$ is smooth in simplices $T_j$ completely contained 
in some sub-domain $\Omega_k$, that is there exists a $k$ with $T_j \subset \Omega_k$. In this case 
we only search for the nearest neighbor stencil in the reduced set $\left\{\mathbf{x}_i | 
\mathbf{x}_i \in \Omega_k\right\}$, but not in the complete set of sampling points 
$\{\mathbf{x}_i\}$. As in the original approach, we ensure that the vertices $\mathbf{x}_{i_j}$ of 
simplex $T_j$ are contained in the nearest neighbor stencil $S_j$. Since $S_j \subset \Omega_k$, we 
can approximate a smooth function by polynomial interpolation with known order of convergence 
$(p_j+1)/d$. Figure \ref{improved_stencil_1} shows the improved stencil for a simplex $T_j$ 
without any kinks inside.

\begin{figure}
	\centering
	\hfill
	\begin{subfigure}{0.4\textwidth}
		\centering
		\includegraphics{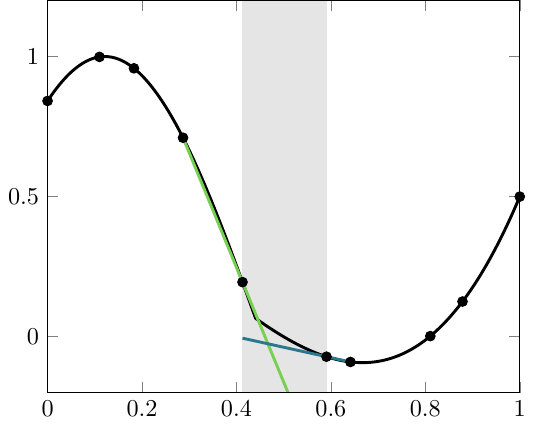}
		\subcaption{Linear approximation.}
	\end{subfigure}
	\hfill
	\begin{subfigure}{0.4\textwidth}
		\centering
		\includegraphics{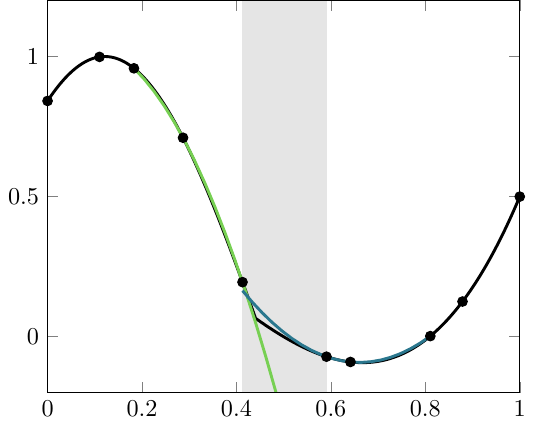}
		\subcaption{Quadratic approximation.}
	\end{subfigure}
	\hfill\mbox{}
	\caption{A linear (a) and quadratic (b) approximation of a kink, each with two stencils. The 
	function $f(x)$ is colored in black, the left hand approximation $g_{j,1}(x)$ in green and the 
	right hand approximation $g_{j,2}(x)$ in blue.}
	\label{approx:kink}
\end{figure}

\paragraph{Case 2.} Suppose simplex $T_j$ is divided by a kink, that is some of its vertices 
$\mathbf{x}_{i_j}$ belong to $\Omega_{j_1}$ and some to $\Omega_{j_2}$, see 
Figure~\ref{improved_stencil_2}. In this case, we search for two nearest neighbor stencils 
$S_{j,1}\subset\Omega_{j_1}$, $S_{j,2}\subset\Omega_{j_2}$ and two approximations 
$g_{j,1}(\mathbf{x})$, $g_{j,2}(\mathbf{x})$, one at each side of the kink. As above we ensure that 
each stencil contains the corresponding vertices $\mathbf{x}_{i_j}$ of $T_j$. Without loss of 
generality, we assume that the kink can be represented for all $\mathbf{x}_i \in S_{j,1} \cup 
S_{j,2}$ as the maximum of both interpolations, i.e.\ 
\begin{align*}
	f(\mathbf{x}_i) = \max\left(g_{j,1}(\mathbf{x}_i),\; g_{j,2}(\mathbf{x}_i)\right).
\end{align*}
Then we extrapolate $g_{j,1}(\mathbf{x})$ and $g_{j,2}(\mathbf{x})$ to simplex $T_j$ and 
approximate $f(\mathbf{x})$ for all $\mathbf{x} \in T_j$ by taking the maximum of both 
approximations
\begin{align*}
	f(\mathbf{x}) \approx g_j(\mathbf{x}) := \max\left(g_{j,1}(\mathbf{x}),\; 
	g_{j,2}(\mathbf{x})\right)
\end{align*}
whereby we obtain an approximation to the kink. Figure \ref{approx:kink} shows a linear and a 
quadratic approximation to a kink in simplex $T_j$. On both stencils $S_{j_1}$ and $S_{j_2}$ the 
function $f$ is smooth. Both approximations $g_{j,1}(\mathbf{x})$ and $g_{j,2}(\mathbf{x})$ converge 
with an order of $(p_j+1)/d$, respectively. Hence, the approximation $g_j(\mathbf{x})$ converges 
with the same order. Note that this holds also true if $g_{j,1}(\mathbf{x})$ and 
$g_{j,2}(\mathbf{x})$ do not intersect in $T_j$. Even if there was not any kink in the function, 
this procedure of computing two approximations and taking the maximum would not affect the 
convergence. This is important because in our application an activated regulator, may cause a kink 
in the flux in some pipes but not in all. \\

Observe that with this new approach we do not need to fulfill the LEC condition (\ref{LEC}) 
anymore. Since we approximate only smooth functions, no oscillations caused by jumps (Gibb's 
phenomenon) will arise. Any oscillations due to Runge's phenomenon will result in a larger error 
estimator and thus in a finer discretization. Indeed, using the LEC limiter would reduce the 
convergence rate if there are some small oscillations in $f(\mathbf{x})$.

\subsection{Refinement Strategies}
While it is possible to construct an approximation for a given set of sampling points, we want to 
start with an initial set of sampling points consisting of the corners and the center of $\Omega$. 
To adaptively refine the discretization  we then successively add new points at those simplices for 
which a to be defined error estimator is the largest. In the end we aim for less points in regions 
were $f(\mathbf{x})$ is flat and more points in regions where $f(\mathbf{x})$ varies more. For an 
adaptive refinement we need on the one hand a strategy of how to add new points and on the other 
hand a reliable error estimator.

\subsubsection{Adding a New Sampling Point}
\begin{figure}[h]
\begin{center}
	\begin{subfigure}{0.4\textwidth}
		\centering
		\includegraphics{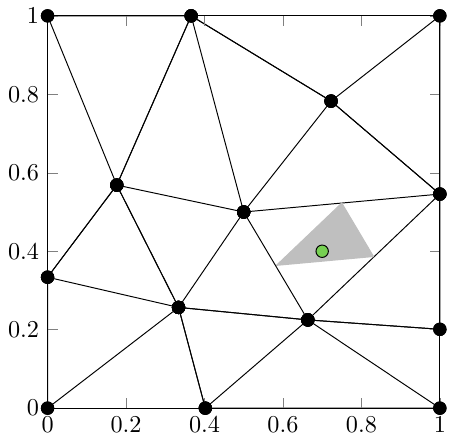}
		\subcaption{New sampling point (green) in subsimplex $T_{\text{sub}_j}$ (light gray).}
	\end{subfigure}
	\hspace{1cm}
	\begin{subfigure}{0.4\textwidth}
		\centering
		\includegraphics{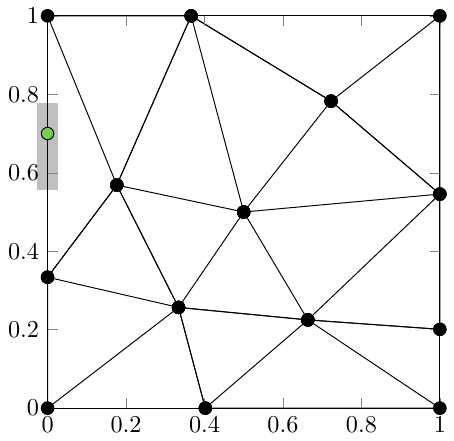}
		\subcaption{New sampling point (green) on the longest edge at the boundary.}
	\end{subfigure}
\end{center}
	\caption{Shown are different refinement strategies for simplices with all facets lying in the 
	interior of the domain $\Omega$ (a) and for simplices with a facet lying at the boundary 
	$\partial \Omega$ (b).}
	\label{refinement}
\end{figure}

In \cite{witteveen2012_2} simplex $T_j$ is refined by sampling a new random point in a subsimplex 
$T_{\text{sub}_j}$. The vertices $\mathbf{x}_{\text{sub}_{j,l}}$ are defined as the centers of the 
faces of simplex $T_j$
\begin{align*}
	\mathbf{x}_{\text{sub}_{j,l}} = \frac{1}{d} \sum_{\genfrac{}{}{0pt}{1}{l^*=0}{l^*\neq l}}^d 
	\mathbf{x}_{i_{j,l^*}}.
\end{align*}
See Appendix \ref{app:a} for an efficient way to sample random points from a uniform distribution 
over some simplex. Figure \ref{refinement}a shows the subsimplex $T_{\text{sub}_j}$ of simplex 
$T_j$. This sampling strategy results in long and flat simplices at the boundary because the new 
sampling point will almost surely not be added at the boundary. Therefore, we use this strategy only 
for simplices without a facet at the boundary. Simplices with a facet at the boundary are refined by 
adding a new sampling point on the middle third of the longest edge, as introduced in 
\cite{witteveen2012_1}. Let $\mathbf{x}_{i_{j,0}}$ and $\mathbf{x}_{i_{j,1}}$ be the endpoints of 
the longest edge of simplex $T_j$, then we define the new sampling point 
$\mathbf{x}_{i_\text{new}}$ as
\begin{align*}
	\mathbf{x}_{i_\text{new}} = \mathbf{x}_{i_{j,0}} + 
	\textstyle\frac{1+u}{3}\;(\mathbf{x}_{i_{j,1}}-\mathbf{x}_{i_{j,0}}),
\end{align*}
where $u$ is a uniformly distributed random variable in $[0,1]$. Figure \ref{refinement}b shows the 
sampling area on the longest edge of a boundary simplex $T_j$.

\subsubsection{Error Estimation}
To refine the simplex with the largest error we need an error estimator since we cannot compute the 
exact error. First, we introduce two newly developed solution-based error estimators and then a 
third already existing error estimator that does not directly depend on the solution. The third one 
is very useful when the function $f(\mathbf{x})$ has more than one output.

\paragraph{Error Estimation Based on a Single Point.} 
In \cite{witteveen2012_1} a solution-based error estimator $\varepsilon_j$ is proposed where the 
square of the hierarchical error $\epsilon_{i_{\text{new},j}} = |f(\mathbf{x}_{i_{\text{new},j}}) 
-g_j(\mathbf{x}_{i_{\text{new},j}})|$ between approximation and function at the new sampling point 
$\mathbf{x}_{i_{\text{new},j}}$ is weighted with the volume of the simplex 
\begin{align}
	\varepsilon_j = \operatorname{vol}(T_j) \cdot \epsilon_{i_{\text{new},j}}^2.
\end{align}
This error estimator has the disadvantage that we need to evaluate the function $f$ at point 
$\mathbf{x}_{i_{\text{new},j}}$, although the point might not be added to the discretization. 
To avoid these useless function evaluations we modify the original error estimator. 
We do not use the hierarchical error in the new sampling point $\mathbf{x}_{i_{\text{new},j}}$, but 
instead in the last added sampling point of simplex $T_j$ before adding it. Let $i_{j^*} = \max_l 
i_{j,l}$ be the index of this last added sampling point and $T_{\text{ref},j^*}$ the simplex which 
was refined by adding $\mathbf{x}_{i_{j^*}}$. Then, the hierarchical error is given by 
\begin{align*}
	\epsilon_{i_{j^*}} = |f(\mathbf{x}_{i_{j^*}}) - g_{\text{ref},j^*} (\mathbf{x}_{i_{j^*}})|
\end{align*}
and we obtain the error estimator
\begin{align*}
	\tilde{\varepsilon}_j & = \operatorname{vol}(T_j) \cdot \epsilon_{i_{j^*}}^2. 
\end{align*} 
By summing up the error estimators for all simplices $\{T_j\}$ we can approximate the root mean 
square error in $\Omega$ by
\begin{align*}
	\tilde{\varepsilon}_\text{rms} = \sqrt{\sum_{j=1}^m \tilde{\varepsilon}_j} = \sqrt{\sum_{j=1}^m 
\operatorname{vol}(T_j) \cdot \epsilon_{i_{j^*}}^2 }.
\end{align*}

\paragraph{Error Estimation Based on Monte Carlo Integration.}
Because we do not want to rely on the error in one single point, we develop a new error estimator.
It is an approximation of the $L_1$ error between $g_j(\mathbf{x})$ and the function $f(\mathbf{x})$ 
in a given simplex $T_j$. For this, we approximate $\varepsilon_j = \|f-g_j\|_{L_1(T_j)}$ by Monte 
Carlo integration, i.e.
\begin{align}
	\varepsilon_j \approx \operatorname{vol}(T_j) \sum_{i=1}^{n_\text{MC}} 
	\frac{|f(\mathbf{x}_{\text{MC},i}) - g_j(\mathbf{x}_{\text{MC},i})|}{n_\text{MC}}
	\label{eps:mc}
\end{align}
at $n_\text{MC}$ randomly drawn Monte Carlo points $\mathbf{x}_{\text{MC},i}$. It is not feasible to 
evaluate $f$ at all $n_{\text{MC}}$ Monte Carlo points because each function evaluation can be an 
expensive simulation. Thus we approximate the right hand side of (\ref{eps:mc}) with the polynomial 
interpolation $\overline{g}_j$ in stencil $S_j$ of degree $p_j-1$. We define
\begin{align*}
	\hat{\varepsilon}_j = \operatorname{vol}(T_j) \sum_{i=1}^{n_\text{MC}} 
	\frac{|g_j(\mathbf{x}_{\text{MC},i}) - \overline{g}_j(\mathbf{x}_{\text{MC},i})| 
	^{(p_j+1)/p_j}}{n_\text{MC}}.
\end{align*}
The exponent $(p_j+1)/p_j$ is necessary since the approximation with $\overline{g}_j$ only leads to 
an order of convergence of $p_j/d$, whereas the approximation $g_j$ converges with order 
$(p_j+1)/d$. Thereby we ensure that the error estimator decreases with the same rate as the true 
error. If $p_j=1$, we define the constant function $\overline{g}_j$ as $\overline{g}_j(\mathbf{x}) 
= \min_{i_j} f(\mathbf{x}_{i_j})$. To obtain an overall error estimation, we sum up the error 
estimators for all simplices $\{T_j\}$
\begin{align*}
	\hat{\varepsilon}_{l_1} = \sum_{j=1}^m \hat{\varepsilon}_j.
\end{align*}

\paragraph{Error Estimation Based on the Theoretical Order of Convergence.}
If one has a function $f(\mathbf{x})$ with a multidimensional output, a solution-based error 
estimator could not be used because one usually does not know how the error scales over different 
outputs. Therefore we use a solution-independent error estimator, as described in 
\cite{witteveen2012_1}. For this consider the definition of the order of convergence
\begin{align*}
	\mathcal{O} = \frac{\log(\varepsilon_0/\varepsilon_j)} 
								{\log(\operatorname{vol}(\Omega)/\operatorname{vol}(\Omega_j))}
\end{align*}
for some reference error $\varepsilon_0$. Then the error $\varepsilon_j$ in simplex $T_j$ is 
proportional to
\begin{align*}
	\varepsilon_j 
	\sim \operatorname{vol}(T_j)^\mathcal{O} 
	= \operatorname{vol}(T_j)^{(p_j+1)/d}.
\end{align*}
Weighting this again with the volume of simplex $T_j$ yields the error estimator
\begin{align*}
 \overline{\varepsilon}_j = \operatorname{vol}(T_j) \cdot 
\varepsilon_j = \operatorname{vol}(T_j)^{(p_j+1)/d+1}.
\end{align*}
It only depends on the volume of simplex $T_j$ and the theoretical order of convergence 
$\mathcal{O}=(p_j+1)/d$. For an overall error estimator we sum again over all simplices
\begin{align*}
	\overline{\varepsilon}_\mathcal{O} = \sum_{j=1}^m \overline{\varepsilon}_j.
\end{align*}

\subsection{Numerical Results for Test Functions}
\label{sec:num_test}
Here, we provide numerical results for smooth and non-smooth functions. To verify the convergence 
rates we calculate the approximation error as the $l_1$ norm between $f(\mathbf{x})$ and 
$g(\mathbf{x})$ evaluated at $n_\text{MC}=10^6$ uniformly in $\Omega$ distributed Monte Carlo 
points:
\begin{align*}
	\varepsilon_{l_1} & = \sum_{i=1}^{n_\text{MC}} \frac{|f(\mathbf{x}_{\text{MC},i}) 
												 - g(\mathbf{x}_{\text{MC},i})|}{n_\text{MC}}.
\end{align*}

\subsubsection{Smooth Functions}
\begin{figure}
	\centering
	\begin{subfigure}{0.33\textwidth}
		\centering
		\includegraphics{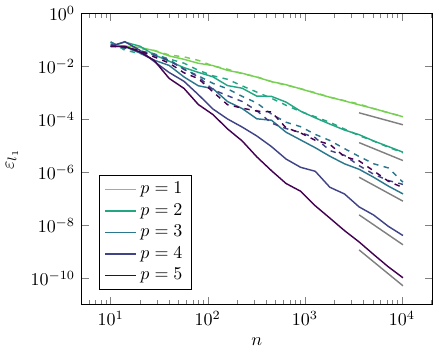}
		\subcaption{$d=2$}
	\end{subfigure}
	\hfill
	\begin{subfigure}{0.33\textwidth}
		\centering
		\includegraphics{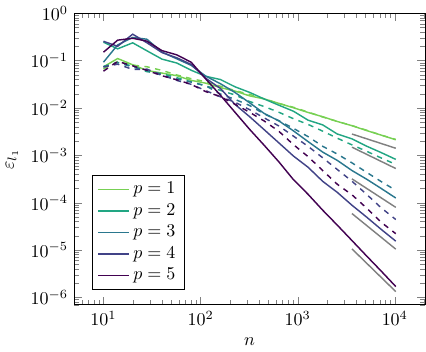}
		\subcaption{$d=3$}
	\end{subfigure}
	\hfill
	\begin{subfigure}{0.33\textwidth}
		\centering
		\includegraphics{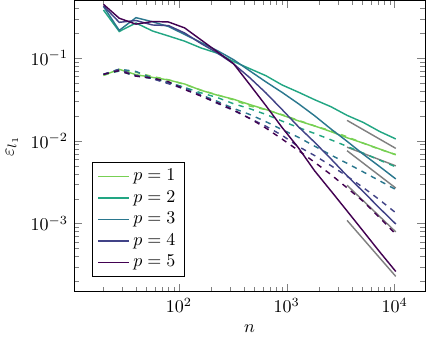}
		\subcaption{$d=4$}
	\end{subfigure}
	\caption{Shown is the $l_1$ error evaluated at $10^6$ random points versus the number $n$ of 
	sampling points for the smooth test function in $d=2,3,4$ dimensions with (dashed lines) and 
	without (solid lines) the LEC condition. The theoretical convergence rates are colored in 
	gray.}
	\label{convergence_smooth_functions}
\end{figure}

First we evaluate the simplex stochastic collocation algorithm with some smooth function $f \in 
C^\infty([0,1]^d)$
\begin{align*}
	f(\mathbf{x}) = \prod_{i=1}^d \sin(\pi x^{(i)}),
\end{align*}
for the Monte Carlo based error estimator $\hat{\varepsilon}_j$ with and without the local extremum 
conserving condition. In Figure \ref{convergence_smooth_functions} we see that the algorithm without 
the LEC condition yields slightly better results for $d \leq 3$. Since the function is smooth, 
oscillations due to kinks or jumps cannot occur. Enforcing the local extremum conservation decreases 
the polynomial degree $p_j$ if the function $f(\mathbf{x})$ itself has some small oscillations in 
simplex $T_j$. This reduction of the polynomial degree is not necessary and impairs convergence. 

\begin{figure}
	\centering
	\begin{subfigure}{0.33\textwidth}
		\centering
		\includegraphics{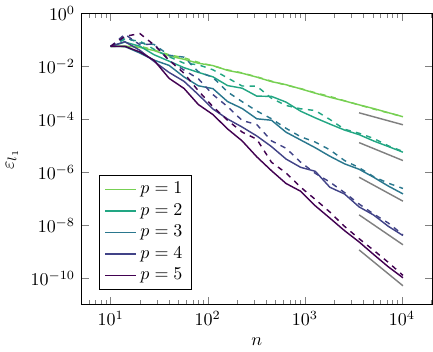}
		\subcaption{$d=2$}
	\end{subfigure}
	\hfill
	\begin{subfigure}{0.33\textwidth}
		\centering
		\includegraphics{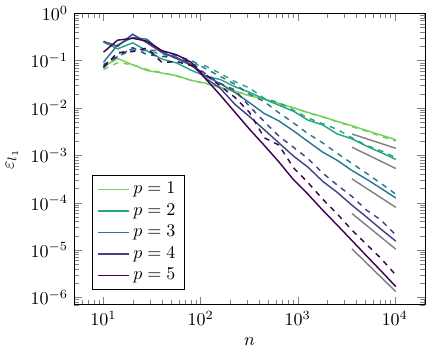}
		\subcaption{$d=3$}
	\end{subfigure}
	\hfill
	\begin{subfigure}{0.33\textwidth}
		\centering
		\includegraphics{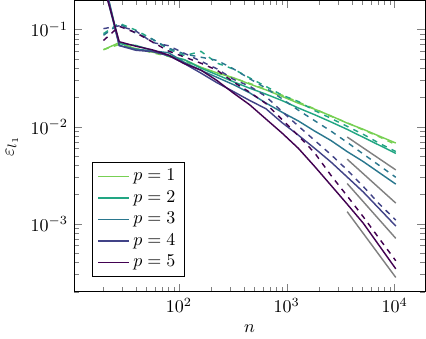}
		\subcaption{$d=4$}
	\end{subfigure}
	\caption{Shown is the $l_1$ error evaluated at $10^6$ random points versus the number $n$ of 
	sampling points for the smooth test function in $d=2,3,4$ dimensions with (dashed lines) and 
	without (solid lines) the assumption of a kink at $f(\mathbf{x})=0.7$. The theoretical 
	convergence rates are colored in gray. Dimensions $d=2$ and $d=3$ are without the LEC condition 
	and $d=4$ is with the $\delta$-LEC condition.}
	\label{convergence_smooth_functions_with_kinks}
\end{figure}

But, with increasing dimension we benefit from using the condition of local extremum conservation in 
the pre-asymptotic behavior. Therefore, we will use a weaker formulation of the local extremum 
conserving condition for dimensions $d\geq 4$ in the following. We will only reduce the polynomial 
degree of the approximation by one, if it does not hold that
\begin{align*}
	\min_{\mathbf{x} \in T_j} g_j(\mathbf{x})+ \delta \geq \min_{\mathbf{x}_i \in T_j} 
	f(\mathbf{x}_i)
	\quad \wedge \quad
	\max_{\mathbf{x} \in T_j} g_j(\mathbf{x})-\delta \leq \max_{\mathbf{x}_i \in T_j} 
f(\mathbf{x}_i),
\end{align*}
with $\delta=0.5(\max_{\mathbf{x}_i \in T_j} f(\mathbf{x}_i)-\min_{\mathbf{x}_i \in T_j} 
f(\mathbf{x}_i))$. This $\delta$-local extremum conserving ($\delta$-LEC) condition allows small 
oscillations in the approximation and improves the pre-asymptotic behavior without affecting the 
convergence. See Figure \ref{convergence_smooth_functions_with_kinks}c for the error in $d=4$ 
dimensions with this weaker condition. 

Above we stated that calculating two approximations on both sides of an assumed kink does not affect 
the convergence rate if in fact there is no kink. In order to verify this, we took the same test 
function $f(\mathbf{x})=\prod_{i=1}^d \sin(\pi x_i)$ and assumed a kink at $f(\mathbf{x}) = 0.7$. 
That is, we check if the function value $f(\mathbf{x}_i)$ is smaller or equal to 0.7, so we can 
assign each sampling point $x_i$ either to $\Omega_1 = \{\mathbf{x}\in\Omega: f(\mathbf{x})<0.7\}$ 
or to $\Omega_2 = \{\mathbf{x}\in\Omega: f(\mathbf{x})=0.7\}$. This replicates the effect of a 
regulator in a gas network. See Figure \ref{convergence_smooth_functions_with_kinks} for the 
results. Assuming a kink yields slightly larger errors, but in all cases the desired convergence 
rates are attained. The difference between assuming and not assuming a kink decreases with 
increasing number $n$ of sampling points. Note that for $d=4$ dimensions we have already used the 
$\delta$-local extremum conservation.

\subsubsection{Non-Smooth Functions}
Consider the test function 
\begin{align*}
	f(\mathbf{x}) = \min\left(\prod_{i=1}^d \sin(\pi x^{(i)}), 0.7\right).
\end{align*}

\begin{figure}
	\centering
	\begin{subfigure}{0.32\textwidth}
		\centering
		\includegraphics{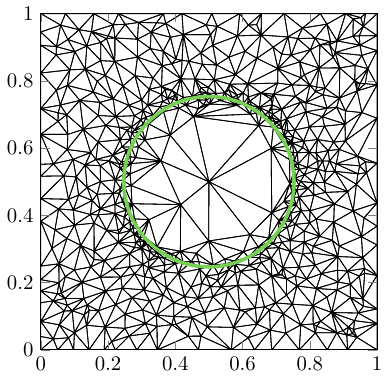}
		\subcaption{$p_j=2$}
	\end{subfigure}
	\hfill
		\begin{subfigure}{0.32\textwidth}
		\centering
		\includegraphics{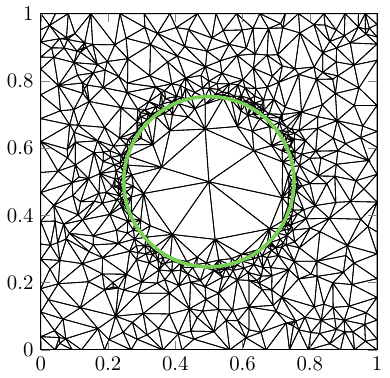}
		\subcaption{$p_j=3$}
	\end{subfigure}
	\hfill
		\begin{subfigure}{0.32\textwidth}
		\centering
		\includegraphics{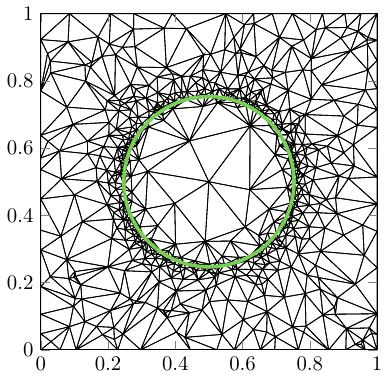}
		\subcaption{$p_j=4$}
	\end{subfigure}
	\caption{Original SSC: Shown is an adaptively refined Delaunay triangulation of $n=640$ sampling 
	points for the non-smooth test function with different polynomial degrees of $p_j=2,3,4$ and the 
	$l_1$ error estimator $\tilde{\varepsilon}_j$. The location of the kink is marked in green.}
	\label{tria_original}
\end{figure}

First we show numerical results for the original simplex stochastic collocation version 
\cite{witteveen2012_2} with the local extremum conservation and no special approximation for kinks.
As expected, enforcing the local extremum conservation reduces the polynomial degree near the kink, 
which results in a larger error estimator and thus in a finer discretization, see Figure 
\ref{tria_original}. The higher the polynomial degree is, the more points are added near the kink. 
We expected this behavior because the smooth part of $f$ can be better approximated with polynomials 
of higher degree, whereas increasing the degree of the interpolating polynomials does not benefit 
approximating the kink. 

\begin{figure}
	\centering
	\begin{subfigure}{0.32\textwidth}
		\centering
		\includegraphics{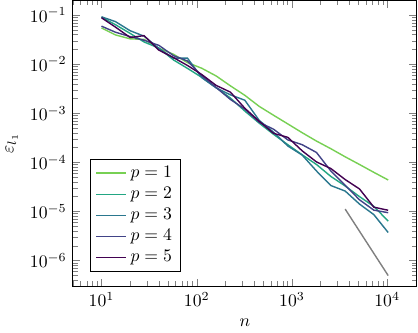}
		\subcaption{$d=2$}
	\end{subfigure}
	\hfill
	\begin{subfigure}{0.32\textwidth}
		\centering
		\includegraphics{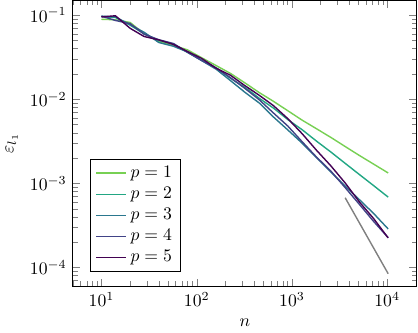}
		\subcaption{$d=3$}
	\end{subfigure}
	\hfill
	\begin{subfigure}{0.32\textwidth}
		\centering
		\includegraphics{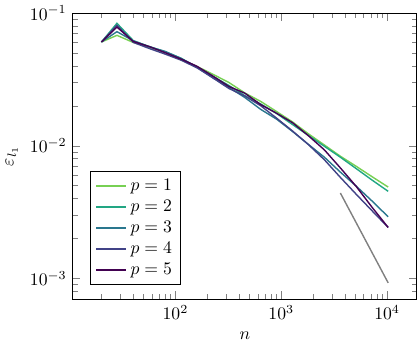}
		\subcaption{$d=4$}
	\end{subfigure}
	\caption{Original SSC: Shown is the $l_1$ error evaluated at $10^6$ random points versus the 
	number $n$ of sampling points for different dimensions $d=2,3,4$ with the $l_1$ error estimator 
	$\tilde{\varepsilon}_j$. The desired convergence rates for a polynomial degree of $p=5$ are 
	plotted in gray.}
	\label{error_original}
\end{figure}

See Figure \ref{error_original} for the convergence rates of the original simplex stochastic 
collocation with the original error estimator $\varepsilon_j$. In $d=2$ dimensions, the desired 
convergence rates are attained for small polynomial degrees $p=1$ and $p=2$. Increasing the 
polynomial degrees up to $p=3$, $p=4$ and $p=5$ does not improve the convergence rate and hence the 
theoretical orders of 2, 2.5 and 3 are not attained. In $d=3$ dimensions the errors for $p=3$, $p=4$ 
and $p=5$ are nearly the same with a maximal order of 1.3 instead of 2. For dimensions larger or 
equal to $d=4$ using polynomials of higher degree is not beneficial and the maximally attained 
order of convergence is 0.75. Therefore, the original simplex stochastic collocation is useless for 
computing statistics of the solution in $d\geq4$ dimensions. For these cases Monte Carlo methods 
provide comparable results with less computational effort. \\

\begin{figure}
	\centering
	\begin{subfigure}{0.32\textwidth}
		\centering
		\includegraphics{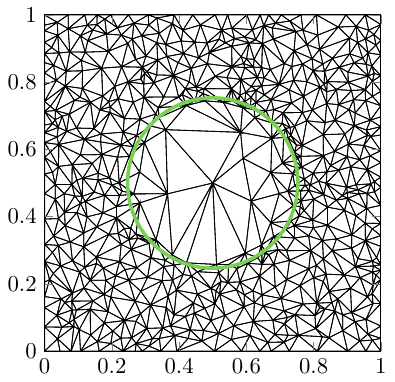}
		\subcaption{Error estimator $\tilde{\varepsilon}_j$}
	\end{subfigure}
	\hfill
	\begin{subfigure}{0.32\textwidth}
		\centering
		\includegraphics{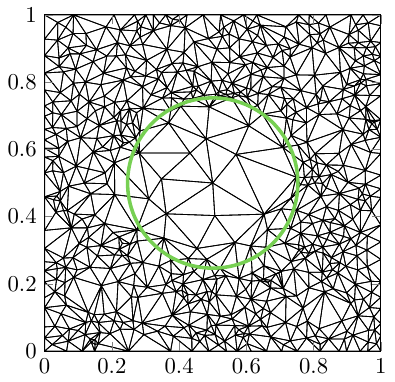}
		\subcaption{Error estimator $\hat{\varepsilon}_j$}
	\end{subfigure}
	\hfill
	\begin{subfigure}{0.32\textwidth}
		\centering
		\includegraphics{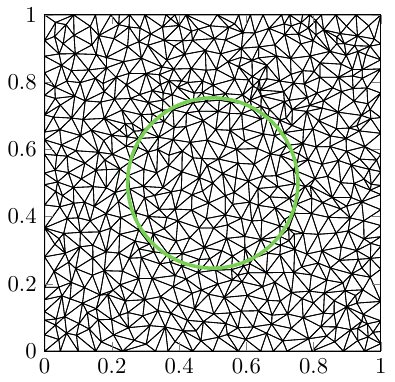}
		\subcaption{Error estimator $\overline{\varepsilon}_j$}
	\end{subfigure}
	\caption{Modified SSC: Shown is an adaptively refined Delaunay triangulation with $n=640$ sampling 
	points for the function $f(\mathbf{x}) = \textstyle \min(\prod_{i=1}^d \sin(\pi x^{(i)}), 0.7)$ in 
	2d with a polynomial degree of $p_j=5$ for the $l_1$ error estimator $\tilde{\varepsilon}_j$ (a), 
	for the root mean square error estimator $\hat{\varepsilon}_j$ (b), and for the 
	function-independent error estimator $\overline{\varepsilon}_j$ (c). The location of the kink is 
	marked in green.}
	\label{trias}
\end{figure}

Now we analyze the modified simplex stochastic collocation method. As before, we check if the 
function value is smaller or equal to 0.7 to simulate a regulator in a gas network. An adaptively 
refined Delaunay triangulation obtained with the $l_1$ error estimator $\hat{\varepsilon}_j$ for 
$p_j=5$ and $n=640$ sampling points can be found in Figure \ref{trias}a. As expected, the sampling 
points are more or less uniformly distributed over the parameter space $\Omega$ where the function 
value is not constant. In the center of our domain where the function value is constant, the areas 
of the triangles are significantly larger. The triangulation in Figure \ref{trias}b, using the root 
mean square error estimator $\tilde{\varepsilon}_j$, looks quite similar: there are fewer triangles 
in the center than around it where the triangles are less uniformly sized as for the estimator 
$\hat{\varepsilon}_j$. In contrast, the triangulation resulting from the function-independent error 
estimator $\overline{\varepsilon}_j$ is uniform. It is not possible to recognize the location of the 
kink, see Figure \ref{trias}c. 

\begin{figure}
	\centering
	\begin{subfigure}{0.325\textwidth}
		\centering
		\includegraphics{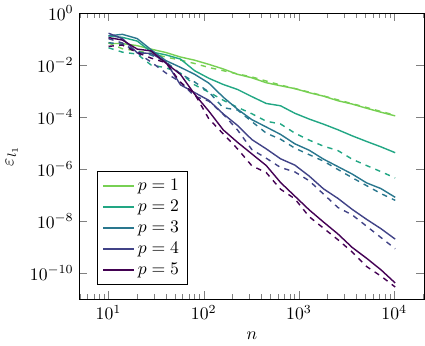}
		\subcaption{$l_1$ error for $\tilde{\varepsilon}_j$, 2d}
	\end{subfigure}
	\hfill
	\begin{subfigure}{0.325\textwidth}
		\centering
		\includegraphics{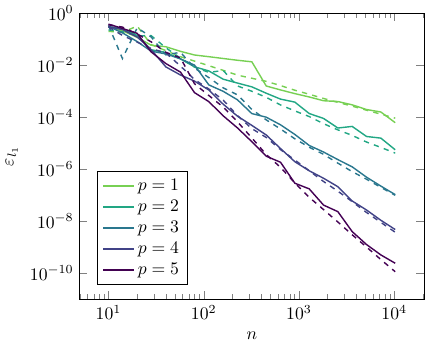}
		\subcaption{$l_1$ error for $\hat{\varepsilon}_j$, 2d}
	\end{subfigure}
	\hfill
	\begin{subfigure}{0.325\textwidth}
		\centering
		\includegraphics{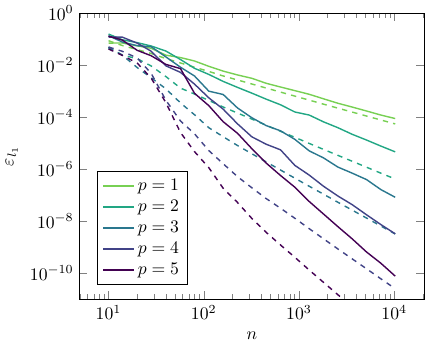}
		\subcaption{$l_1$ error for $\overline{\varepsilon}_j$, 2d}
	\end{subfigure}\\[2em]
	\begin{subfigure}{0.325\textwidth}
		\centering
		\includegraphics{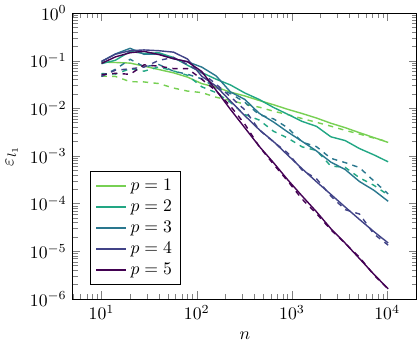}
		\subcaption{$l_1$ error for $\tilde{\varepsilon}_j$, 3d}
	\end{subfigure}
	\hfill
	\begin{subfigure}{0.33\textwidth}
		\centering
		\includegraphics{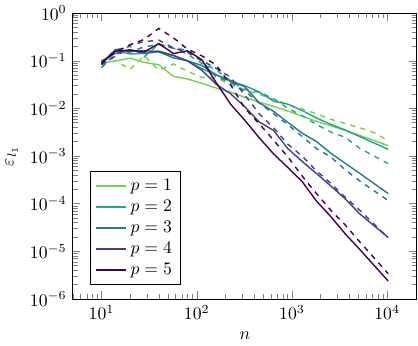}
		\subcaption{$l_1$ error for $\hat{\varepsilon}_j$, 3d}
	\end{subfigure}
	\hfill
	\begin{subfigure}{0.325\textwidth}
		\centering
		\includegraphics{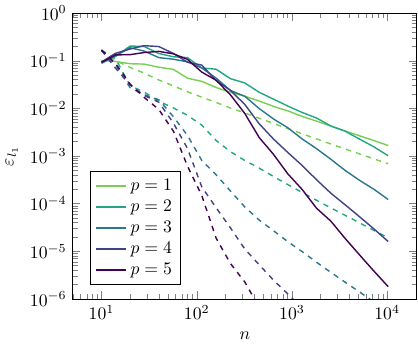}
		\subcaption{$l_1$ error for $\overline{\varepsilon}_j$, 3d}
	\end{subfigure}\\[2em]
	\begin{subfigure}{0.325\textwidth}
		\centering
		\includegraphics{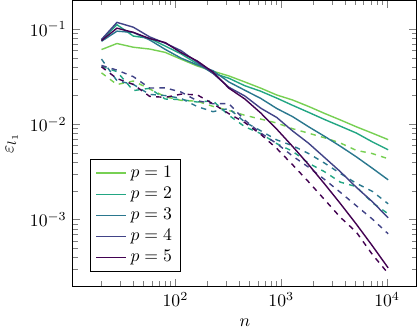}
		\subcaption{$l_1$ error for $\tilde{\varepsilon}_j$, 4d}
	\end{subfigure}
	\hfill
	\begin{subfigure}{0.325\textwidth}
		\centering
		\includegraphics{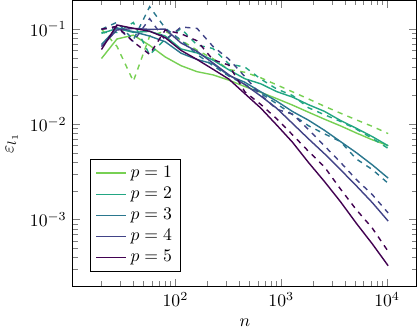}
		\subcaption{$l_1$ error for $\hat{\varepsilon}_j$, 4d}
	\end{subfigure}
	\hfill
	\begin{subfigure}{0.325\textwidth}
		\centering
		\includegraphics{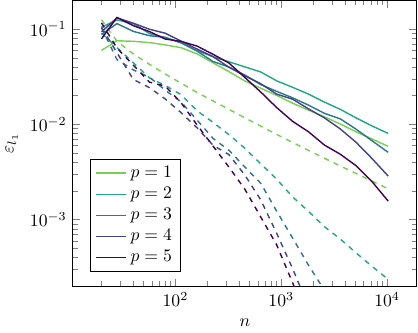}
		\subcaption{$l_1$ error for $\overline{\varepsilon}_j$, 4d}
	\end{subfigure}
	\caption{Modified SSC: Shown is the $l_1$ error evaluated at $10^6$ random points (solid) and the 
	error estimator (dashed) versus the number $n$ of interpolation points for the $l_1$ error 
	estimator $\tilde{\varepsilon}_j$ (a), the root mean square error estimator $\hat{\varepsilon}_j$ 
	(b), and the error estimator $\overline{\varepsilon}_j$ (c).}
	\label{error_estimators}
\end{figure}

In all shown dimensions $d=2,3,4$ nearly all theoretical convergence rates of $\varepsilon_{l_1}$ 
evaluated at $n_\text{MC}=10^6$ Monte Carlo points are attained for the $l_1$ error estimator as 
well as for the root mean square error estimator $\hat{\varepsilon}_j$, and the error estimator 
$\overline{\varepsilon}_j$, cf.\ Figure \ref{error_estimators}. The $l_1$ error estimator yields 
the best results and the smoothest convergence. The pointwise error estimator 
$\tilde{\varepsilon}_j$ yields comparable results and both error estimators can be used as a 
reliable stopping criterion. The total errors reached with error estimator 
$\overline{\varepsilon}_j$ in $d=2$ and $d=3$ dimensions look quite similar, but as expected the 
estimated overall error differs greatly from the real error because it is not solution-based. So 
this error estimator should not be used as stopping criterion. Moreover, this is also the reason for 
the worse results in $d=4$ dimensions. The error estimator $\overline{\varepsilon}_j$ overestimates 
the real error in simplices where the polynomial degree has been reduced for fulfilling the 
$\delta$-LEC condition. Omitting the $\delta$-LEC condition in this case would decrease the error 
for a large number of sampling points but at the expense of a worse pre-asymptotic. Comparing these 
total errors with those obtained with the original simplex stochastic collocation method and the 
original pointwise error estimator $\varepsilon_j$, shows that the modification yields significantly 
better results. The total error for the maximal number of points was improved from $3\cdot10^{-6}$ 
to $4\cdot10^{-11}$ in two dimensions, from $2\cdot10^{-4}$ to $2\cdot10^{-6}$ in three dimensions, 
and from $2\cdot 10^{-3}$ to $3\cdot10^{-4}$ in four dimensions.

\subsubsection{Multiple Refinements}
\begin{figure}
	\centering
	\begin{subfigure}{0.32\textwidth}
		\centering
		\includegraphics{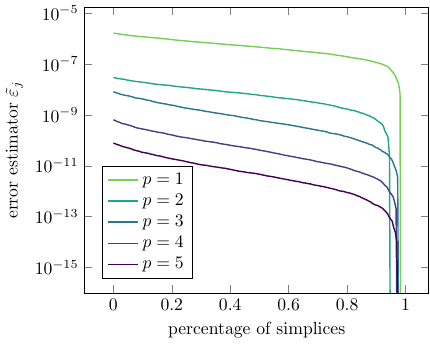}
		\subcaption{d=2, error estimator $\tilde{\varepsilon}_j$}
	\end{subfigure}
	\hfill
	\begin{subfigure}{0.32\textwidth}
		\centering
		\includegraphics{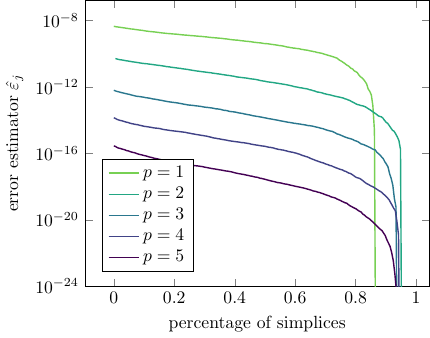}
		\subcaption{d=2, error estimator $\hat{\varepsilon}_j$}
	\end{subfigure}
	\hfill
	\begin{subfigure}{0.32\textwidth}
		\centering
		\includegraphics{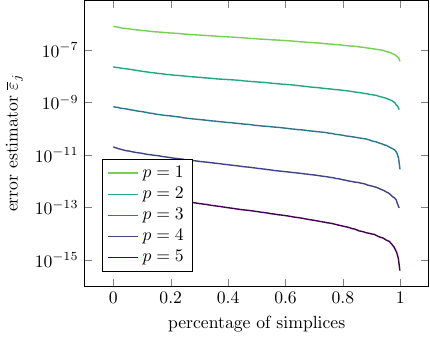}
		\subcaption{d=2, error estimator $\overline{\varepsilon}_j$}
	\end{subfigure}\\[2em]	
	
	\begin{subfigure}{0.32\textwidth}
		\centering
		\includegraphics{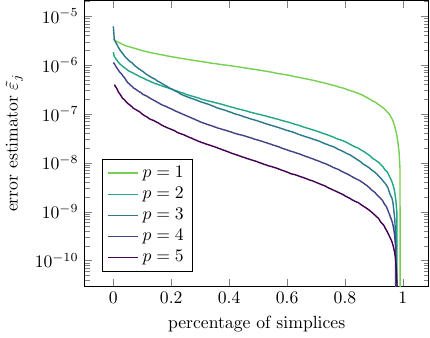}
		\subcaption{d=3, error estimator $\tilde{\varepsilon}_j$}
	\end{subfigure}
	\hfill
	\begin{subfigure}{0.32\textwidth}
		\centering
		\includegraphics{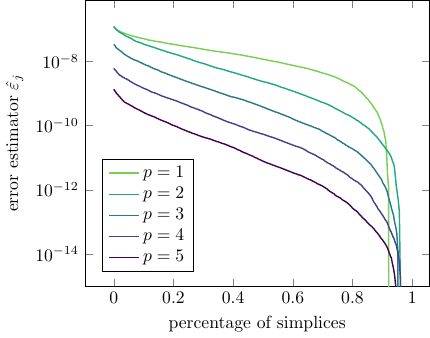}
		\subcaption{d=3, error estimator $\hat{\varepsilon}_j$}
	\end{subfigure}
	\hfill
	\begin{subfigure}{0.32\textwidth}
		\centering
		\includegraphics{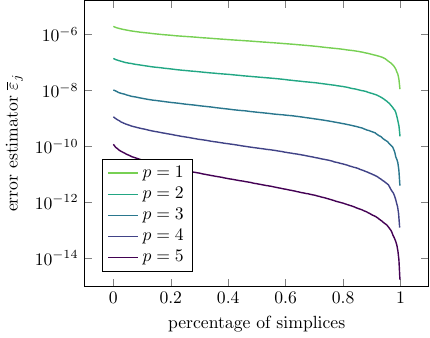}
		\subcaption{d=3, error estimator $\overline{\varepsilon}_j$}
	\end{subfigure}\\[2em]
		\centering
	\begin{subfigure}{0.32\textwidth}
		\centering
		\includegraphics{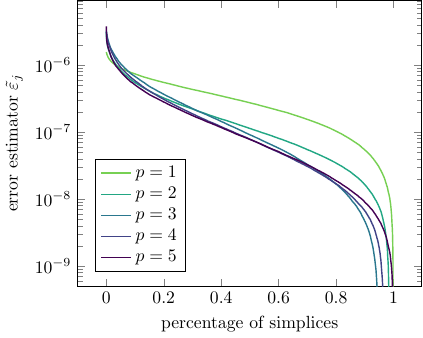}
		\subcaption{d=4, error estimator $\tilde{\varepsilon}_j$}
	\end{subfigure}
	\hfill
	\begin{subfigure}{0.32\textwidth}
		\centering
		\includegraphics{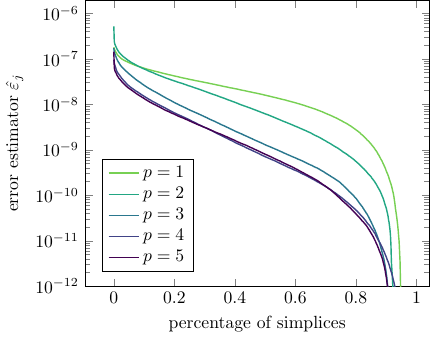}
		\subcaption{d=4, error estimator $\hat{\varepsilon}_j$}
	\end{subfigure}
	\hfill
	\begin{subfigure}{0.32\textwidth}
		\centering
		\includegraphics{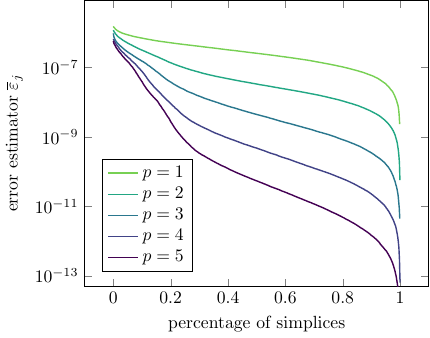}
		\subcaption{d=4, error estimator $\overline{\varepsilon}_j$}
	\end{subfigure}
	\caption{Shown is the distribution of the $l_1$ error estimator over the percentage of simplices 
	for different error estimators. In all cases, $m=1000$ sampling points are used for the 
	triangulation.}
	\label{error_distribution}
\end{figure}

In order to parallelize the refinement, at each step the $m_\text{ref}\geq1$ simplices with the 
largest error estimator can be refined. Thus, the function evaluations for the new sampling points 
can be done simultaneously and the expensive update of the Delaunay triangulation needs just to be 
done once instead of $m_\text{ref}$ times. Figure \ref{error_distribution} shows the three error 
estimators versus the percentage of simplices. The results are comparable for all dimensions and all 
error estimators, the only exception is error estimator $\overline{\varepsilon}_j$ in $d=4$ 
dimensions. In this case, the estimated error in 20 \% of the simplices is higher then expected. 
These are exactly the simplices where the polynomial degree was reduced for fulfilling the 
$\delta$-LEC condition. The all other cases, the error estimator slowly decreases over most 
simplices independent of dimension, polynomial degree, and type of error estimator. Only for a small 
percentage of simplices the error estimator is significantly smaller than for the rest. Thereby it 
is reasonable to add several sampling points at once.

\begin{figure}
	\centering
	\begin{subfigure}{0.32\textwidth}
		\centering
		\includegraphics{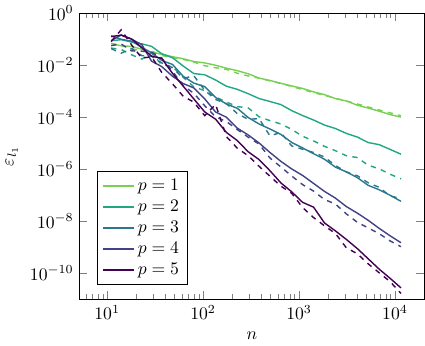}
		\subcaption{d=2, $m_\text{ref}=0.3n$}
	\end{subfigure}
	\hfill
	\begin{subfigure}{0.32\textwidth}
		\centering
		\includegraphics{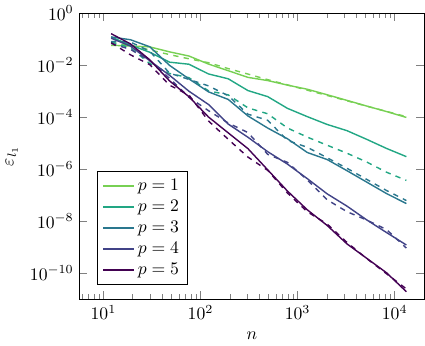}
		\subcaption{d=2, $m_\text{ref}=0.6n$}
	\end{subfigure}
	\hfill
	\begin{subfigure}{0.32\textwidth}
		\centering
		\includegraphics{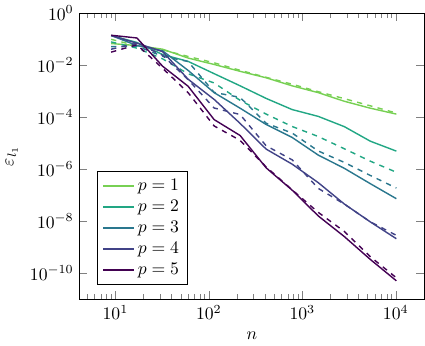}
		\subcaption{d=2, $m_\text{ref}=0.9n$}
	\end{subfigure}\\[2em]
	
	\begin{subfigure}{0.32\textwidth}
		\centering
		\includegraphics{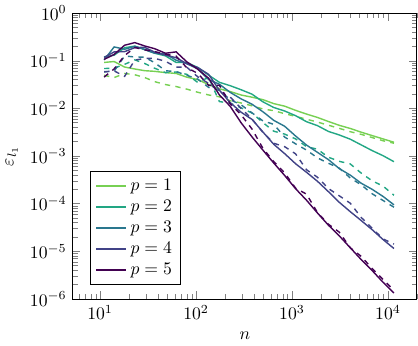}
		\subcaption{d=3, $m_\text{ref}=0.3n$}
	\end{subfigure}
	\hfill
	\begin{subfigure}{0.32\textwidth}
		\centering
		\includegraphics{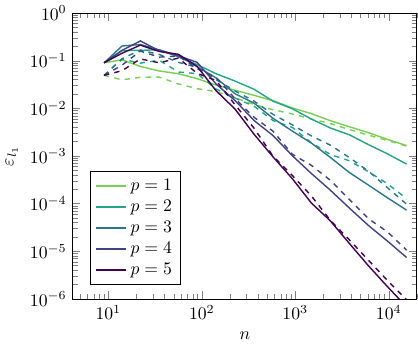}
		\subcaption{d=3, $m_\text{ref}=0.6n$}
	\end{subfigure}
	\hfill
	\begin{subfigure}{0.32\textwidth}
		\centering
		\includegraphics{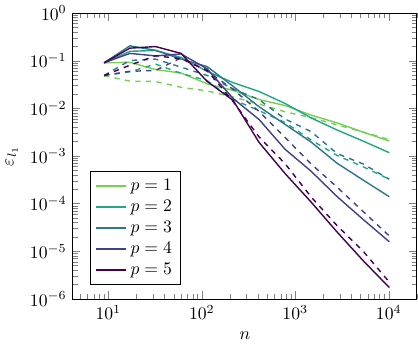}
		\subcaption{d=3, $m_\text{ref}=0.9n$}
	\end{subfigure}\\[2em]
	
	\begin{subfigure}{0.32\textwidth}
		\centering
		\includegraphics{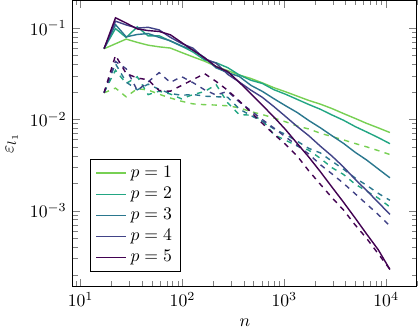}
		\subcaption{d=4, $m_\text{ref}=0.3n$}
	\end{subfigure}
	\hfill
	\begin{subfigure}{0.32\textwidth}
		\centering
		\includegraphics{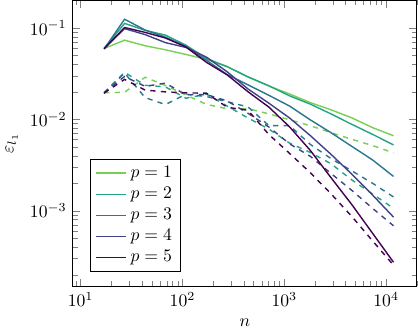}
		\subcaption{d=4, $m_\text{ref}=0.6n$}
	\end{subfigure}
	\hfill
	\begin{subfigure}{0.32\textwidth}
		\centering
		\includegraphics{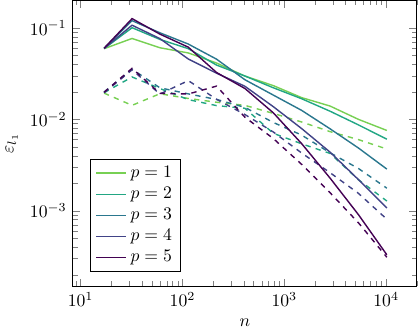}
		\subcaption{d=4, $m_\text{ref}=0.9n$}
	\end{subfigure}

	\caption{Shown is the $l_1$ error evaluated at $10^6$ random points (solid) and the error 
	estimator (dashed) versus interpolation points $n$ for $l_1$ error estimator 
	$\tilde{\varepsilon}_j$. At each refinement step 30\% (left column), 60\% (central column), and 
	90\% (right column) of the old sampling points were added.}
	\label{error_multiple}
\end{figure}

Figure \ref{error_multiple} shows the convergence rates for multiple refinements where we used the 
Monte Carlo based error estimator $\tilde{\varepsilon}_j$ and at each step added $0.3n$, $0.6n$, or 
$0.9n$ points, respectively, to the current discretization consisting of $n$ sampling points. Since 
the number of newly added sampling points does not influence the convergence, it is reasonable to 
refine multiple simplices to save computational time. 

\clearpage
\subsection{Statistics of the Approximated Function}
When simulating gas networks some input data can be uncertain like the pressure of the injected gas 
at input nodes or the flux of the extracted gas at demand nodes. The response of the gas network to 
these uncertainties are expressed by the pressure and temperature at nodes and the flux through 
pipes. We are interested in statistics of these physical quantities like the expected value, 
variance, or median. The cumulative density function (cdf) can be used to determine the probability 
that a production-related critical value, for example the maximum pressure a pipe can withstand, 
will be exceeded.

\paragraph{Expectation.} The expectation of a function $f(\mathbf{x})$ of a random variable 
$\mathbf{x}\in\Omega$ with the density function $\rho(\mathbf{x})$ can be approximated by using the 
approximations $g_j(\mathbf{x})$ on the simplices $T_j$ in the evaluation of a quadrature rule $Q$ 
for the approximation $\tilde{f}$ of $f$, i.e.\
\begin{align*}
	\mathbb{E}[f] = \int_\Omega f(\mathbf{x}) \rho(\mathbf{x}) \:\textrm{d}\mathbf{x} \approx Q(\tilde{f}) = \sum_{j=1}^m Q(g_j).
\end{align*}
The quadrature rule $Q$ can be a Monte Carlo integration or a Gaussian quadrature, where we do not 
evaluate $f$ for any quadrature point, but only the approximations $g_j$, which are cheap to 
evaluate.

\paragraph{Variance.} The variance of a function $f(\mathbf{x})$ of a random variable 
$\mathbf{x}\in\Omega$ with the density function $\rho(\mathbf{x})$ is defined as the squared 
distance of the function from its mean. We approximate the variance in the same way as the 
expectation, namely by using Monte Carlo integration or Gaussian quadrature to calculate the 
integrals, i.e.
\begin{align*}
	\mathbb{V}[f] = \int_\Omega \left(f(\mathbf{x})-\mathbb{E}[f]\right)^2 \rho(\mathbf{x}) 
	\textrm{d}\mathbf{x} = \mathbb{E}[f^2] - \mathbb{E}[f]^2 \approx Q(\tilde{f}^2) - Q(\tilde{f})^2. 
\end{align*}

\paragraph{Convergence.} The absolute error of the expectation $|\mathbb{E}[f] - Q(\tilde{f})|$ can 
be estimated as
\begin{align*}
	\left|\mathbb{E}[f] - Q(\tilde{f})\right| 
	& \leq \left| \mathbb{E}[f] - \mathbb{E}[\tilde{f}] \right|
		+ \left| \mathbb{E}[\tilde{f}] - Q(\tilde{f}) \right| \\
	& \leq \int_\Omega \left| f(\mathbf{x}) - \tilde{f}(\mathbf{x})\right| 
		\rho(\mathbf{x}) \:\d \mathbf{x} + \left| \mathbb{E}[\tilde{f}] - Q(\tilde{f}) \right| \\
	& \leq \int_\Omega \underbrace{\left| f(\mathbf{x})-\tilde{f}(\mathbf{x}) 
		\right|_\infty}_{\varepsilon_I(f)} \rho(\mathbf{x}) \:\d \mathbf{x}  + \underbrace{\left| 
		\mathbb{E}[\tilde{f}] - Q(\tilde{f}) \right|}_{\varepsilon_Q(f)} \\
	& = \varepsilon_I(f) + \varepsilon_Q(f).
\end{align*}
The interpolation error $\varepsilon_I(f)$ can be estimated by (\ref{eq:convergence}). If we choose 
the quadrature formula such that the quadrature error $\varepsilon_Q(f)$ is at most of the same 
order of magnitude as the interpolation error $\varepsilon_I(f)$, then the approximation 
$Q(\tilde{f})$ of the expected value $\mathbb{E}[f]$ converges also with an order of $(p+1)/d$, 
provided that the partial derivatives are bounded. \\ 

The same rate can be obtained for the variance, if the function and all partial derivatives are 
bounded. Using the triangle inequality we get the following two terms:
\begin{align}
	\left| \mathbb{V}[f] - (Q(\tilde{f}^2) - Q(\tilde{f})^2) \right|  &
	\leq \left| \mathbb{E}[f^2] - Q(\tilde{f}^2) \right| + \left| \mathbb{E}[f]^2 - Q(\tilde{f})^2 
\right|.
	\label{eq:error_variance}
\end{align}
Analogously to the expectation, the first term can be estimated by
\begin{align*}
	\left| \mathbb{E}[f^2] - Q(\tilde{f}^2) \right| \leq \varepsilon_I(f^2) + \varepsilon_Q(f^2).
\end{align*}
With $|f^2-\tilde{f}^2|\leq|f-\tilde{f}|\: |f+\tilde{f}|\leq |f-\tilde{f}|\:(|f|+|\tilde{f}|)$ we 
obtain
\begin{align*}
	\left| \mathbb{E}[f^2] - Q(\tilde{f}^2) \right| \leq (|f|+|\tilde{f}|) \varepsilon_I(f) + 
\varepsilon_Q(f^2).
\end{align*}
Next we consider the second term of (\ref{eq:error_variance}):
\begin{align*}
	\left| \mathbb{E}[f]^2 - Q(\tilde{f})^2 \right|
	& \leq \left| \mathbb{E}[f] - Q(\tilde{f}) \right| \left(|\mathbb{E}[f]|+|Q(\tilde{f})|\right)\\
	& \leq \left(\varepsilon_I(f) + 
	\varepsilon_Q(f)\right)\left(|\mathbb{E}[f]|+|Q(\tilde{f})|\right).
\end{align*}
Assuming bounded $f, \tilde{f}, \mathbb{E}[f], Q(\tilde{f})\leq C$, we get the following result for 
the error of the variance
\begin{align*}
	\left| \mathbb{V}[f] - (Q(\tilde{f}^2) - Q(\tilde{f})^2) \right|
	&\leq 2C\: \varepsilon_I(f) + \varepsilon_Q(f^2) + 2C\:\varepsilon_I(f) + 2C\:\varepsilon_Q(f) \\
	&\leq 4C\: \varepsilon_I(f) + 2C\:\varepsilon_Q(f) + \varepsilon_Q(f^2).
\end{align*}

Hence, by choosing the quadrature formula such that the quadrature errors $\varepsilon_Q(f)$ and 
$\varepsilon_Q(f^2)$ are at most of the same order of magnitude as the interpolation error 
$\varepsilon_I(f)$, yields again an order of convergence of $(p+1)/d$. 

\paragraph{CDF.} For approximating the cumulative density function $\mathbb{P}[f(\mathbf{x}) 
\leq 
y]$, we discretize the value range of the approximation $g$ with equidistant nodes $y_0, y_1, 
\ldots, y_n$ where $y_i = \min(g) + ih$ and $h=(\max(g)-\min(g))/n$. For each node $y_i$ we 
determine the maximal domain $\Omega_i\subseteq \Omega$ such that $g(\mathbf{x})\leq y_i$ for all 
$\mathbf{x} \in \Omega_i$. With the probabilities of these domains we obtain the function values of 
the cumulative density function because it holds 
\begin{align*}
	\mathbb{P}[g(\mathbf{x})\leq y_i] = \mathbb{P}[\Omega_i].
\end{align*}
As a last step we interpolate the cumulative density function between the nodes, e.g.\ with 
piecewise linear polynomials. Note that the interpolation must be monotonically increasing because 
otherwise the resulting function does not fulfill the requirements of a cumulative density function.

\section{Uncertainty Quantification for Gas Network Simulation}
\label{sec:gas_net_results}
Today, natural gas contributes significantly to many countries energy supply, where it is used to 
provide heat and power. Additionally natural gas is an input for producing plastics and chemicals in 
industry. A large number of scenario analyses are necessary to ensure a secure and reliable 
operation of a gas network. Since usually these scenarios cannot be easily tested, they are replaced 
by simulations. Uncertainties arise in the withdrawn amount of gas of each customer. For example, it 
is then of great interest whether the gas network can meet the demand when all customers need a lot 
of gas at once and how likely a failure is. For this forward propagation we use the method of 
simplex stochastic collocation to approximate and integrate high-dimensional functions.

\subsection{Euler Equations for Pipes}
A gas network is modeled with nodes and edges. The edges represent pipes or other network elements 
such as valves, control valves, heaters, or compressors. Gas flow through a single pipe of length 
$L$ with diameter $A$ is described by the Euler equations, a set of partial differential equations 
\cite{koch2015, lurie2008, schmidt2016}. The first equation is the continuity equation 
\begin{align}
 \partial_t\rho + \partial_x (\rho v) = 0,
 \label{eq:continuity}
\end{align}
following from the conservation of mass. The law of momentum conservation 
\begin{align}
	\frac{1}{A} \partial_t q + \partial_x (\rho v^2) + g \rho \partial_x h + \partial_x p + 
	\frac{\lambda}{2D} \rho |v|v = 0,
	\label{eq:momentum}
\end{align}
specifies the pressure loss along the pipe due to weight, pressure, and frictional forces. The 
equation of state 
\begin{align}
 p = z(p,T) \rho R_s T
 \label{eq:state}
\end{align}
is necessary to describe the state of a real compressible gas for a given set of values for 
temperature $T$, density $\rho$, and pressure $p$. The first law of thermodynamics must be taken 
into account to describe any heat transfer process. A solution to this system of equations can be 
found analytically if we assume a stationary and isothermal gas flow \cite{schmidt2016}. Analogously 
to Kirchhoff's law, the mass must be conserved at junctions where several pipes are connected. At 
supply nodes the incoming gas pressure is given, whereas at demand nodes the extracted mass flow. If 
a gas network consists of pipes only, the solution of the pressure, density, and temperature at 
nodes and the gas flow in pipes is sufficiently smooth. But a real gas network also contains more 
complicated elements. For an overview over other elements and the corresponding equations see 
\cite{fuchs2018}. In that work errors due to model assumptions were investigated, for realistic 
situations these can be in the order of $10^{-4}$. 

\subsection{Kinks due to pressure regulation}
\begin{figure}
	\centering
	\includegraphics{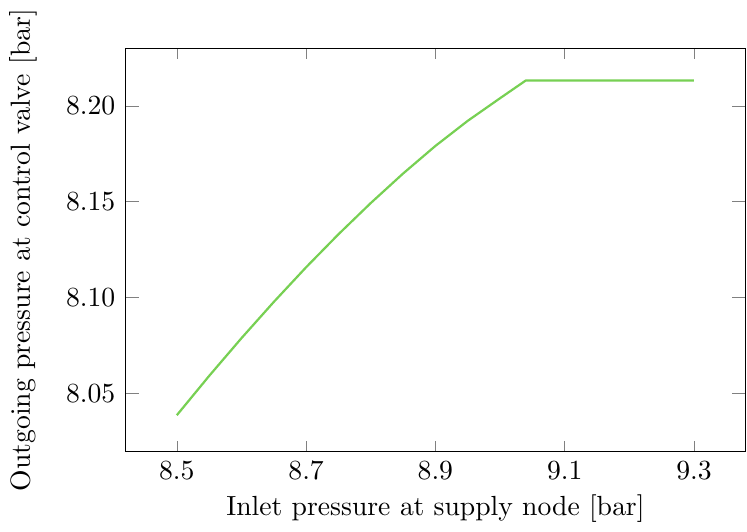}
	\caption{A kink in the solution resulting from pressure regulation.}
	\label{fig:motivation_kink}
\end{figure}

Usually, the pressure in transport pipes is significantly larger than the maximum allowable 
operating pressure in distributional pipes. Due to this reason the network needs pressure control 
valves that adjust the outgoing pressure if the incoming pressure exceeds a preset limit. 
Unfortunately, the more complicated elements impair the smoothness of the solution. For example, a 
pressure control valve causes kinks in the solution. Increasing the pressure at a supply node 
increases the pressure after a control valve until the preset pressure is reached, but afterwards 
the pressure remains constant, see Figure \ref{fig:motivation_kink}. We do not know in advance where 
the kink is located, but after the simulation run we know if a control valve is active or not. This 
information is necessary to use our improved simplex stochastic collocation. 

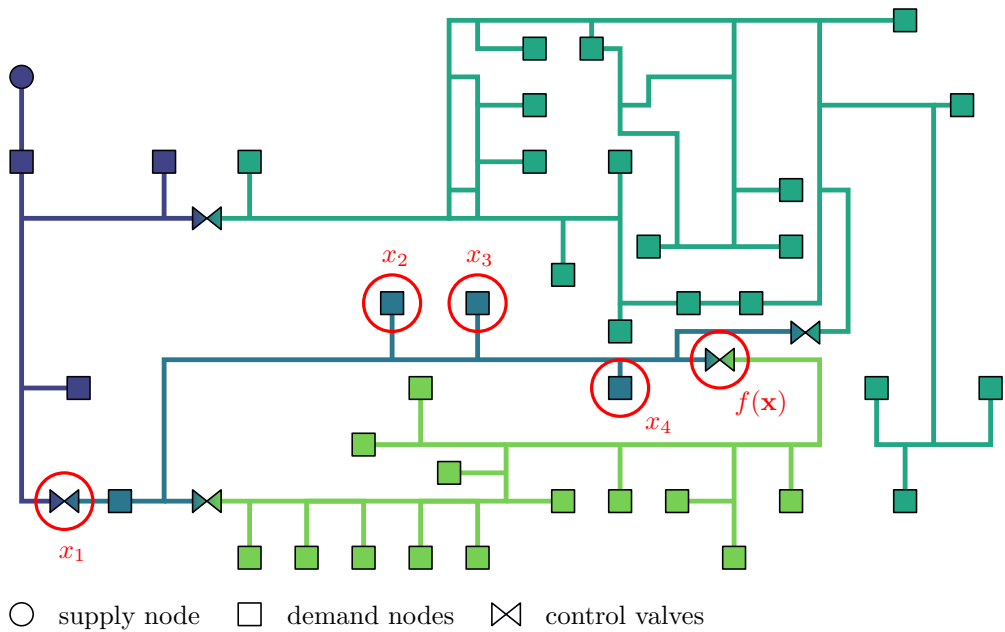
\begin{figure}
	\centering
	\begin{overpic}[scale=0.75]{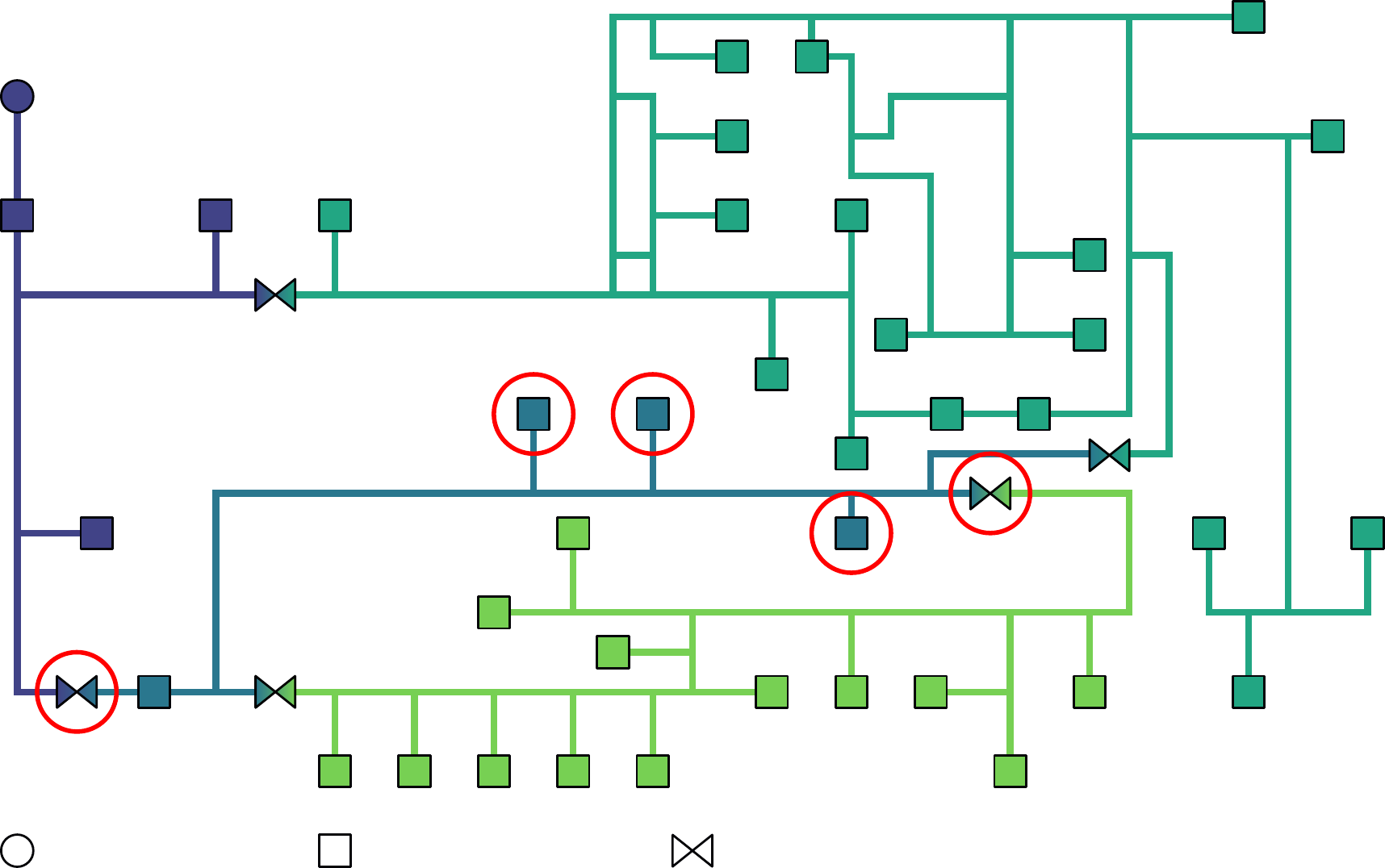} 
	\put(5,0.3){supply node}
	\put(28,0.3){demand nodes}
	\put(54,0.3){control valves}
	\put(5,7){$\textcolor{red}{x_1}$}
	\put(37.5, 37){$\textcolor{red}{x_2}$}
	\put(46, 37){$\textcolor{red}{x_3}$}
	\put(64, 20){$\textcolor{red}{x_4}$}
	\put(73, 22){$\textcolor{red}{f(\mathbf{x})}$}
	\end{overpic}
	\caption{Test gas network with one supply node, 37 demand nodes, and five pressure control 
	valves.}
	\label{fig:network}
\end{figure}

In this section we apply our new version of simplex stochastic collocation to a real gas network 
where we used \cite{clees2016} to simulate the gas flow. The network has one supply node, 37 demand 
nodes, several pipes, and five control valves which reduce the high pressure of about 27 bar at the 
supply node stepwise to pressures of around 16, 8, and 4 bar at the demand nodes. See Figure 
\ref{fig:network} for a schematic drawing of the network. Different pressure levels are shown in 
different colors. In all tests, the quantity of interest is the outgoing pressure $f(\mathbf{x})$ 
at the right control valve. Depending on the uncertain parameters of outgoing pressure $x_1$ at the 
left valve, and the amount of withdrawn gas at demand nodes $x_2$, $x_3$, and $x_4$, the right valve 
is in an active or bypass mode. This is checked by comparing the outgoing pressure with the preset 
pressure. The lower the outgoing pressure $x_1$ is, and the higher the withdrawn amount of gas is, 
the lower is the incoming pressure $f(\mathbf{x})$ at the right control valve.
\vspace{-0.5em}

\subsection{Input Uncertainties in Two Dimensions}
\begin{figure}
	\centering
	\hfill
	\begin{subfigure}{0.45\textwidth}
		\centering
		\includegraphics{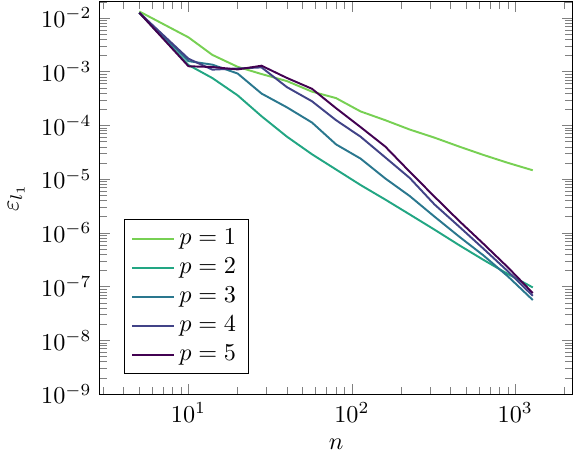}
		\vspace{-0.5em}
		\subcaption{Original SSC}
	\end{subfigure}
	\hfill\hfill
	\begin{subfigure}{0.45\textwidth}
		\centering
		\includegraphics{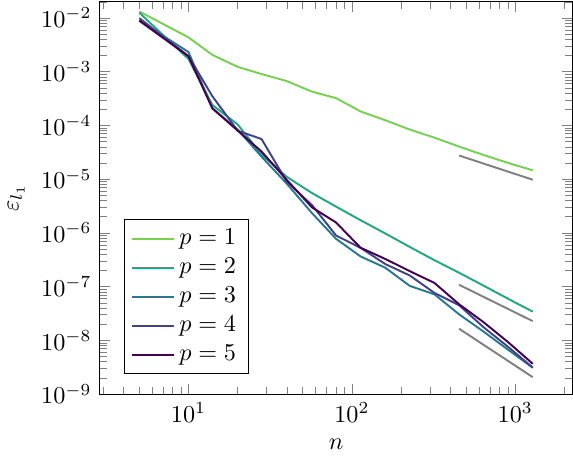}
		\vspace{-0.5em}
		\subcaption{New SSC}
	\end{subfigure}
	\hfill\mbox{}
	\caption{$d=2$. The $l_1$ error evaluated at $10^6$ random points versus the number $n$ of 
	interpolation points with $l_1$ error estimator $\tilde{\varepsilon}_j$ for the original SSC 
	without kink information (a), and for the new version with kink information (b).}
	\label{fig:gas_2d_error}
\end{figure}

\begin{figure}
	\centering
	\hfill
	\begin{subfigure}{0.45\textwidth}
		\centering
		\includegraphics{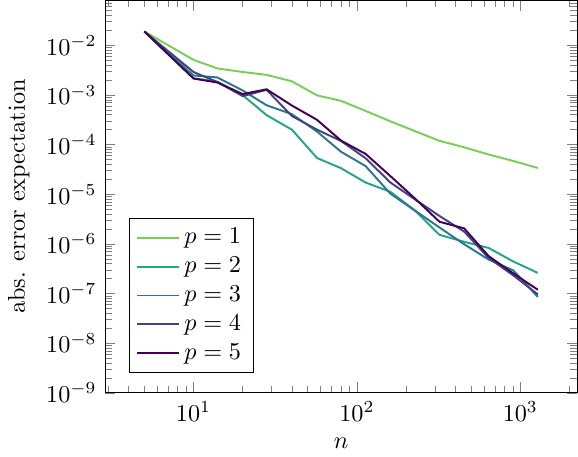}
		\vspace{-0.5em}
		\subcaption{Original SSC}
	\end{subfigure}
	\hfill\hfill
	\begin{subfigure}{0.45\textwidth}
		\centering
		\includegraphics{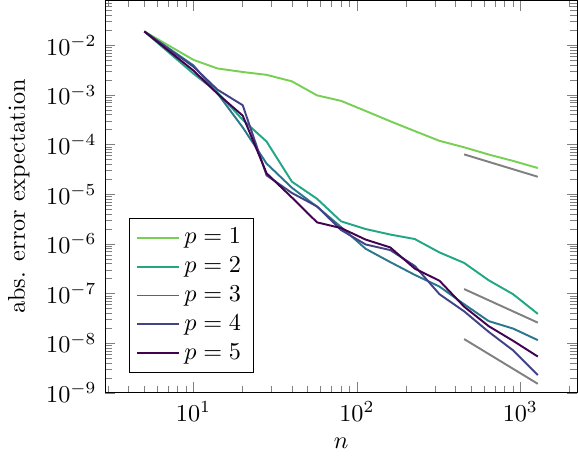}
		\vspace{-0.5em}
		\subcaption{New SSC}
	\end{subfigure}
	\hfill\mbox{}
	\caption{$d=2$. The absolute error in the expected value versus the number $n$ of interpolation 
	points with $l_1$ error estimator $\tilde{\varepsilon}_j$ for the original SSC without kink 
	information (a), and for the new version with kink information (b).}
	\label{fig:gas_2d_expectation}
\end{figure}

First, we vary the outgoing pressure $x_1$ of the left control valve uniformly between 8.5 bar and 
9.5 bar, and the demanded power $x_2$ uniformly between 160 MW and 200 MW. The remaining powers are 
fixed, in particular $x_3=250$ MW and $x_4 = 17$ MW. See Figure \ref{fig:gas_2d_error} for a 
comparison of the results using the original simplex stochastic collocation (a) with those from our 
improved one (b). The new version yields better results than the original one, the smallest error 
reached is two orders of magnitude smaller. In the original version it makes no difference whether 
polynomials of degree $p=2$ or higher are used. In the new version the desired convergence rates 
(marked in gray) for $p=1$, $p=2$, and $p=3$ are obtained. Increasing the polynomial degree to $p=4$ 
or $p=5$ yields no improvement in the rate. Similar results are valid for the expected value, where 
a reference value was computed with a polynomial degree of $p=5$ and $m=5120$ sampling points, see 
Figure \ref{fig:gas_2d_expectation}. The original simplex stochastic collocation needs $m\approx 50$ 
sampling points to achieve an accuracy of $10^{-4}$, about the order of the model error, whereas the 
new version only needs $m\approx 30$ sampling points. 

\begin{figure}
	\hfill
	\begin{subfigure}[b]{0.46\textwidth}
		\centering
		\includegraphics{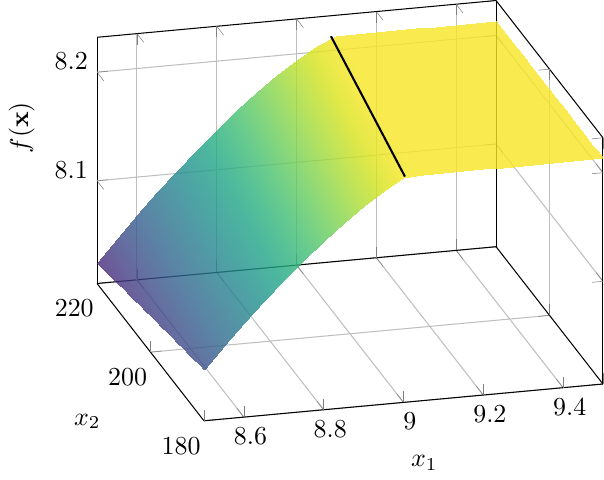}
		\caption{Surface plot}
	\end{subfigure}
	\hfill
	\begin{subfigure}[b]{0.46\textwidth}
	 	\hspace{15mm}
		\begin{overpic}[scale=0.32]{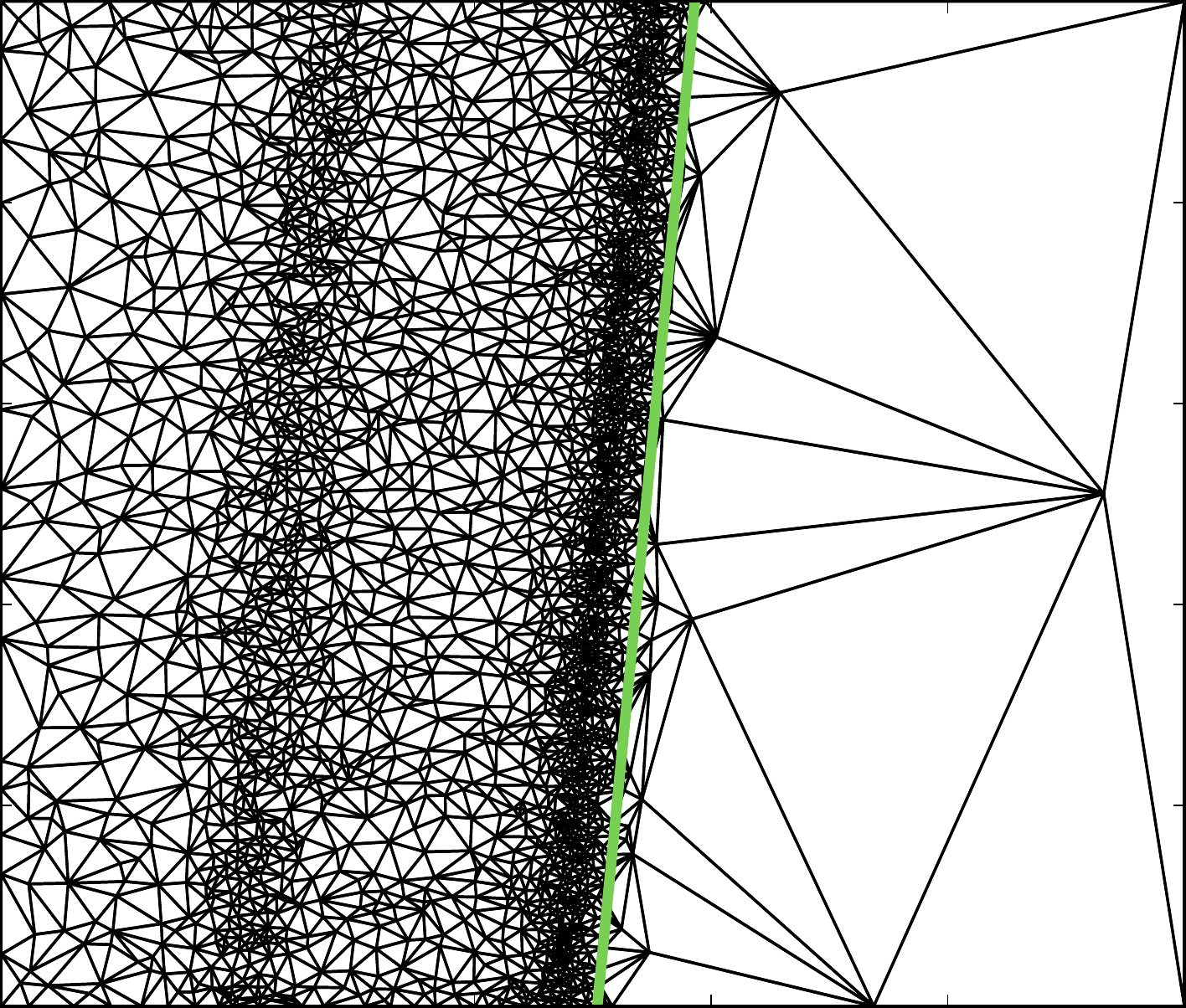}
			\put(-5,-8){\scriptsize 8.5}
			\put(15,-8){\scriptsize 8.7}
			\put(35,-8){\scriptsize 8.9}
			\put(55,-8){\scriptsize 9.1}
			\put(75,-8){\scriptsize 9.3}
			\put(95,-8){\scriptsize 9.5}
			
			\put(-13, -2){\scriptsize 180}
			\put(-13, 20){\scriptsize 190}
			\put(-13, 42){\scriptsize 200}
			\put(-13, 62){\scriptsize 210}
			\put(-13, 83){\scriptsize 220}
		\end{overpic}
		\vspace{8mm}
	 	\caption{Triangulation}
	\end{subfigure}
	\hfill\mbox{}\\
	\hfill
	\begin{subfigure}[b]{0.45\textwidth}
		\centering
		\includegraphics{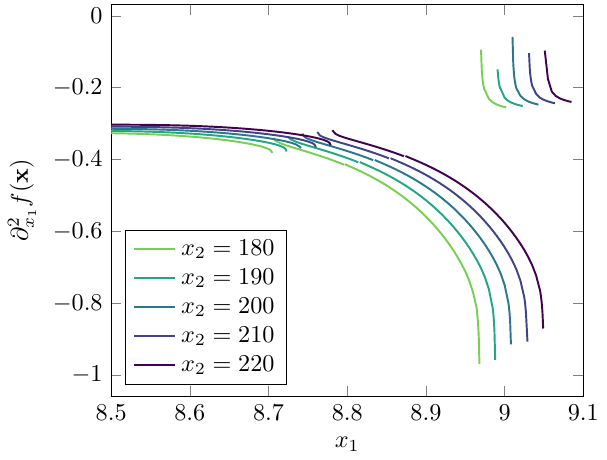}
		\vspace{-1em}
		\caption{Second derivative}
	\end{subfigure}
	\hfill
	\begin{subfigure}[b]{0.45\textwidth}
	 	\centering
		\includegraphics{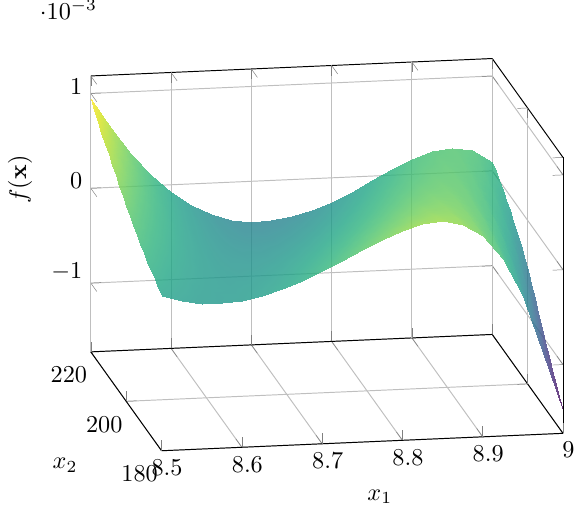}
		\caption{Difference to quadratic regression}
	\end{subfigure}
	\caption{Function $f(\mathbf{x})$ resulting from gas network simulation.}
	\label{fig:function_gas}
\end{figure}
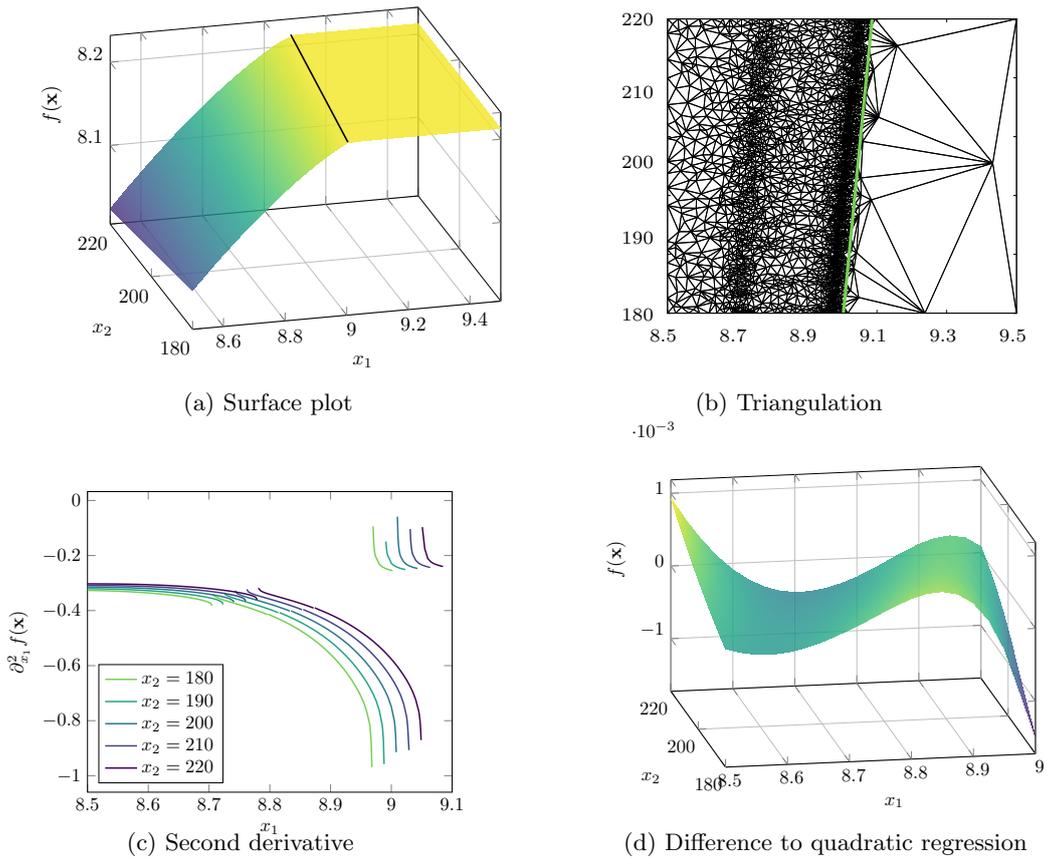

The question is why we cannot improve the convergence rate using higher order polynomials as it was 
the case for the synthetic test function in the previous section? To answer this question see 
Figure \ref{fig:function_gas}. The surface plot of the function (a) is inconspicuous. But if we look 
at the resulting triangulation (b), we can see that there are two areas in which the error estimator 
places an unexpected number of points. This indicates discontinuities. Because of this, we 
approximated the second order partial derivative $\partial_{x_1} f(\mathbf{x})$ with the second 
order finite difference quotient (c) and, hence, we can see two jumps in the second partial 
derivative which explain the poor convergence results. After investigation with the developers of 
the solver MYNTS~\cite{clees2016}, it can be determined that these jumps are not caused by the 
physical properties of gas flow but by the specific numerical treatment of the underlying solver to 
obtain convergence. It is not predictable where these arise, so we do not have a possibility to  
adapt the method. In principal, these jumps due to the numerical treatment in MYNTS could be 
avoided, but this would effect the overall solution process and convergence to a solution would not 
any longer be guaranteed without further measures. Since the employed version of the solver is 
completely sufficient in its accuracy and behavior for current industrial applications, there is so 
far no practical need to overhaul it. 

Because the original version does not have any information about the kink in the function 
$f(\mathbf{x})$, the convergence rates for polynomial degrees of $p\geq 2$ are the same. The new 
version has some information about the kink in the function, but no information about the kink in 
the first derivative (corresponding to the jump in the second derivative) and, therefore, the 
convergence rates for $p\geq 3$ are the same. At first glance, we only improved the order of 
convergence from $1.5$ to $2$, but at second glance, we see that our new simplex stochastic 
collocation has a significantly better pre-asymptotic behavior. This is due to the fact that the 
linear and quadratic terms of $f(\mathbf{x})$ contribute most, whereas the higher order terms are 
only of magnitude $10^{-3}$. See Figure \ref{fig:function_gas}d for the difference between the 
function $f(\mathbf{x})$ and a quadratic regression at the left side of the kink.

\subsection{Input Uncertainties in Three Dimensions}
\begin{figure}
	\centering
	\hfill
	\begin{subfigure}{0.45\textwidth}
		\centering
		\includegraphics{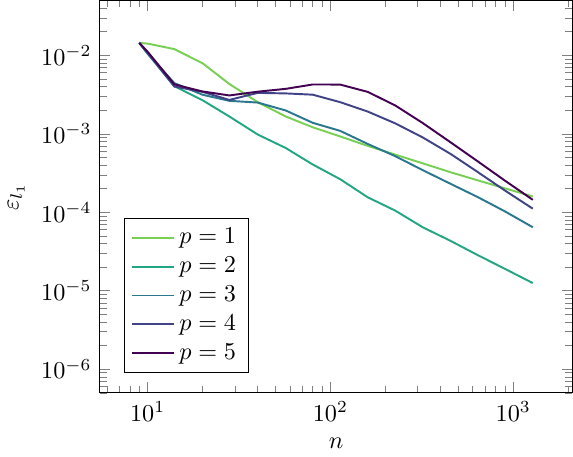}
		\subcaption{Original SSC}
	\end{subfigure}
	\hfill\hfill
	\begin{subfigure}{0.45\textwidth}
		\centering
		\includegraphics{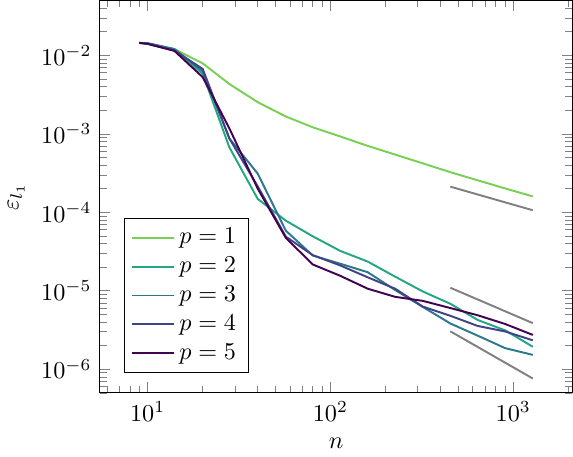}
		\subcaption{New SSC}
	\end{subfigure}
	\hfill\mbox{}
	\caption{$d=3$. The $l_1$ error evaluated at $10^6$ random points versus the number $n$ of 
	interpolation points with $l_1$ error estimator $\tilde{\varepsilon}_j$ for the original SSC 
	without kink information (a), and for the new version with kink information (b).}
	\label{fig:gas_3d_error}
\end{figure}

\begin{figure}
	\centering
	\hfill
	\begin{subfigure}{0.45\textwidth}
		\centering
		\includegraphics{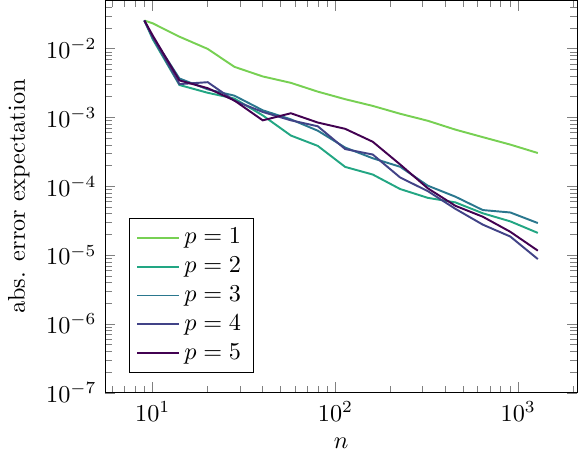}
		\subcaption{Original SSC}
	\end{subfigure}
	\hfill\hfill
	\begin{subfigure}{0.45\textwidth}
		\centering
		\includegraphics{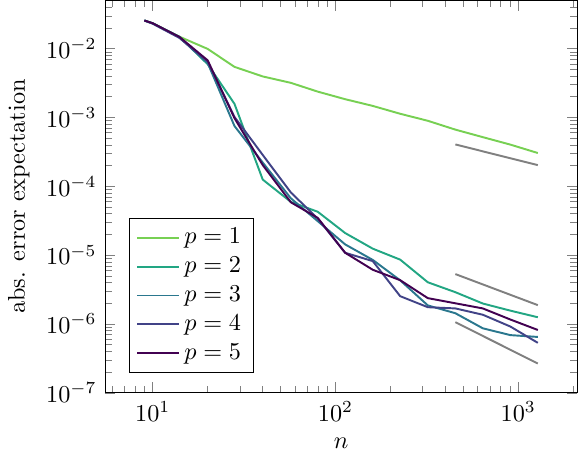}
		\subcaption{New SSC}
	\end{subfigure}
	\hfill\mbox{}
	\caption{$d=3$. The absolute error in the expected value versus the number $n$ of interpolation 
	points with $l_1$ error estimator $\tilde{\varepsilon}_j$ for the original SSC without kink 
	information (a), and for the new version with kink information (b).}
	\label{fig:gas_3d_expectation}
\end{figure}

In addition to the first two uncertain parameters, we now add a third one. The power $x_3$ of the 
withdrawn gas at the marked demand node is uniformly varied between 230 MW and 250 MW. See Figure 
\ref{fig:gas_3d_error}a for the error of the original stochastic simplex collocation. The best 
convergence rate is obtained for a polynomial degree of $p=2$ and increasing the degree results in a 
larger error estimate. This is not the case for our modified simplex stochastic collocation, see 
\ref{fig:gas_3d_error}b. The error is in the same order of magnitude for all polynomial degrees 
$p=2,3,4,5$ and converges with an order of 1. Here, the good pre-asymptotic behavior can be seen 
even better than in $d=2$ dimensions. Figure \ref{fig:gas_3d_expectation} shows the convergence 
results for the expectation. As in $d=2$ dimensions, the reference value is computed with the new 
version of the simplex stochastic collocation, a polynomial degree of five, and $n=5120$ 
interpolation points. For the original version (a), the difference between different polynomial 
degrees is not as large as predicted by the error estimator. The rate is of the same order of 
magnitude as for the new simplex stochastic collocation (b), but the new version benefits from the 
explicit kink approximation in the pre-asymptotic. Hence, only $m\approx50$ instead of $m\approx200$ 
sampling points are necessary to obtain an error of $10^{-4}$. 

\vspace{-0.8em}

\subsection{Input Uncertainties in Four Dimensions}
\begin{figure}[h]
	\centering
	\hfill
	\begin{subfigure}{0.45\textwidth}
		\centering
		\includegraphics{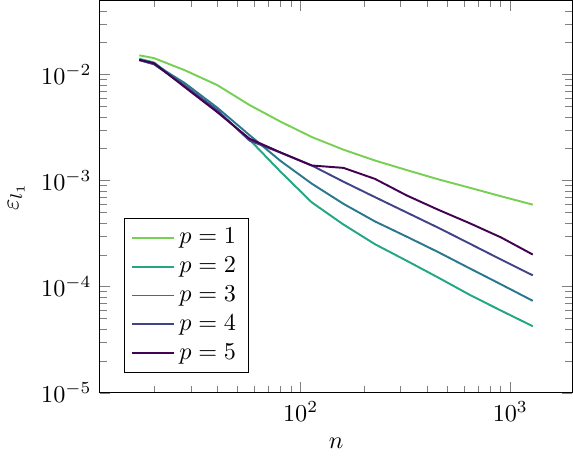}
		\vspace{-0.5em}
		\subcaption{Original SSC}
	\end{subfigure}
	\hfill\hfill
	\begin{subfigure}{0.45\textwidth}
		\centering
		\includegraphics{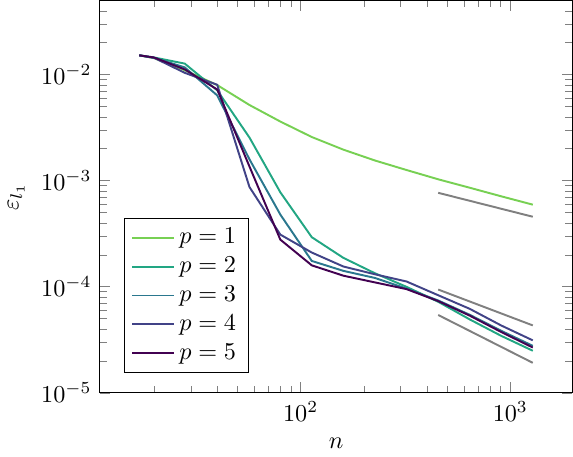}
		\vspace{-0.5em}
		\subcaption{New SSC}
	\end{subfigure}
	\hfill\mbox{}
	\caption{$d=4$. The $l_1$ error evaluated at $10^6$ random points versus the number $n$ of 
	interpolation points with $l_1$ error estimator $\tilde{\varepsilon}_j$ for the original SSC 
	without kink information (a), and for the new version with kink information (b).}
	\label{fig:gas_4d_error}
\end{figure}

\begin{figure}[h]
	\centering
	\hfill
	\begin{subfigure}{0.45\textwidth}
		\centering
		\includegraphics{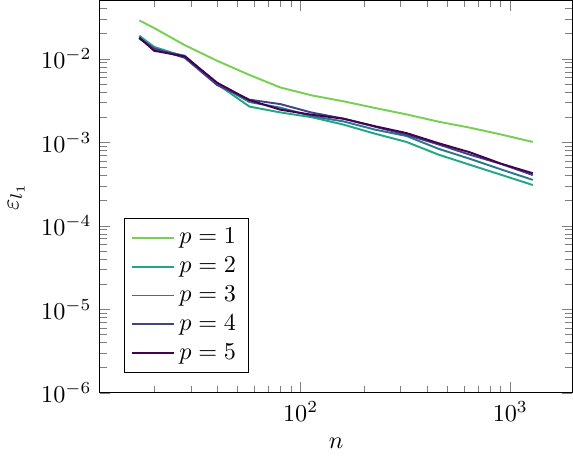}
		\vspace{-0.5em}
		\subcaption{original SSC}
	\end{subfigure}
	\hfill\hfill
	\begin{subfigure}{0.45\textwidth}
		\centering
		\includegraphics{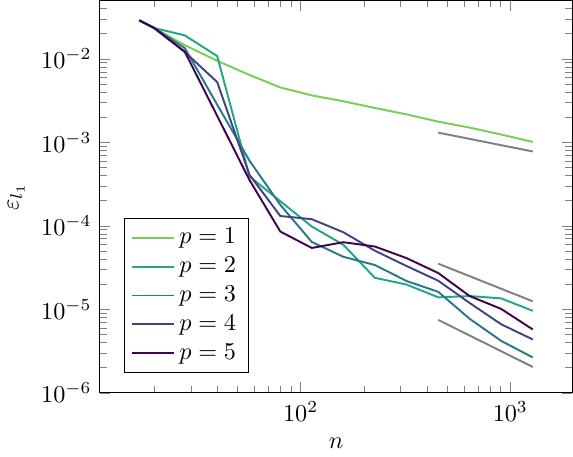}
		\vspace{-0.5em}
		\subcaption{new SSC}
	\end{subfigure}
	\hfill\mbox{}
	\caption{$d=4$. The absolute error in the expected value versus the number $n$ of interpolation 
	points with $l_1$ error estimator $\tilde{\varepsilon}_j$ for the original SSC without kink 
	information (a), and for the new version with kink information (b).}
	\label{fig:gas_4d_expectation}
\end{figure}

Lastly, we add an uncertainty at the power $x_4$ of the withdrawn gas at the third marked demand 
node. The power uniformly varies between 10 MW and 30 MW. As in $d=2$ and $d=3$ dimensions, the 
estimated $l_1$ error of the original simplex stochastic collocation increases with increasing 
polynomial degree, see Figure \ref{fig:gas_4d_error}a. This difference is no longer visible in the 
error of the expectation, where all polynomial degrees result in errors of same order of magnitude, 
see Figure \ref{fig:gas_4d_expectation}a. Again, the new version of the stochastic simplex 
collocation yields better results because of the better pre-asymptotic behavior. There is no visible 
benefit from using polynomials of degree $p\geq3$, but the obtained order of convergence is 1. To 
achieve an error of $10^{-4}$, we only need $m\approx100$ sampling points, whereas the original 
version does not reach this error with $m\approx1000$ sampling points.

\subsection{Comparison to Other Methods}
\begin{figure}[h]
	\centering
	\hspace*{-0.3em}
	\begin{subfigure}[b]{0.33\textwidth}
		\centering
		\includegraphics{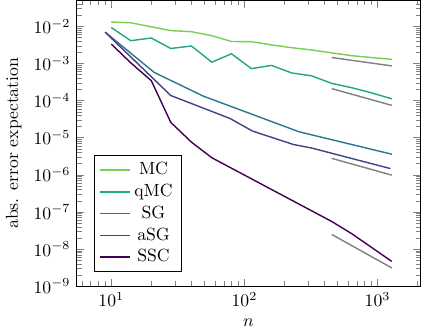}
		\vspace*{-0.5em}
		\subcaption{$d=2$}
	\end{subfigure}
	\hspace{-0.3em}
	\begin{subfigure}[b]{0.33\textwidth}
		\centering
		\includegraphics{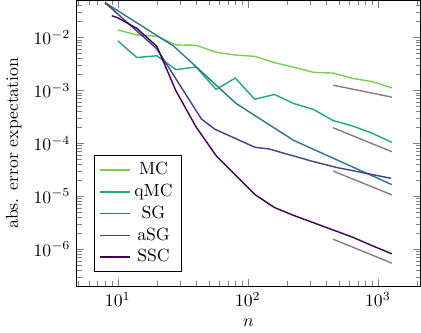}
		\vspace*{-0.5em}
		\subcaption{$d=3$}
	\end{subfigure}
	\hspace{-0.3em}
	\begin{subfigure}[b]{0.33\textwidth}
		\centering
		\includegraphics{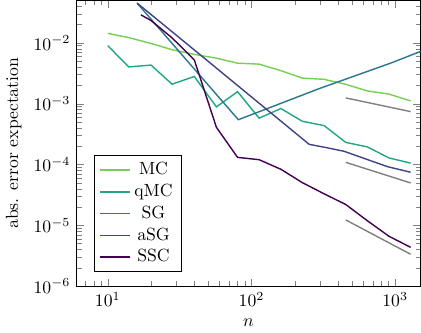}
		\vspace*{-0.5em}
		\subcaption{$d=4$}
	\end{subfigure}
	\hspace*{-0.3em}
	\caption{The absolute error of the expected value versus the number $n$ of interpolation points 
	for Monte Carlo integration, quasi-Monte Carlo integration with the Halton sequence, regular 
	sparse grids, spatially adaptive sparse grids, and simplex stochastic collocation with a 
	polynomial degree of five.}
	\label{fig:gas_comparison}
\end{figure}

Finally, we compare our new simplex stochastic collocation method with other common integration 
methods for computing an expected value. The convergence plots are shown in Figure 
\ref{fig:gas_comparison} for dimensions $d=2$, $d=3$, and $d=4$. The Monte Carlo quadrature does not 
make any requirements on the integrand, therefore, the theoretical order of convergence of 1/2 is 
obtained in all dimensions. The quasi-Monte Carlo quadrature \cite{niederreiter1992} rule with 
Halton points \cite{halton1960, owen2006} yields better results. In $d=2$ and $d=3$ dimensions an 
order of approximately 1 is reached, whereas in $d=4$ dimensions the order is only 3/4. For 
sufficiently smooth integrands, sparse grid quadrature provides even better convergence. Since the 
considered integrand here is only in $\mathcal{C}^{0}(\Omega)$, it is quite interesting how well 
sparse grids perform. We use a regular and a spatially adaptive sparse grid \cite{pflueger2010, 
pflueger2012} with polynomials of degree five. The spatially adaptive variant allows you to place 
more points near singularities or discontinuities \cite{jakeman2011}. In $d=2$ dimensions, both 
sparse grids yield better results than the quasi-Monte Carlo quadrature, but with the same order of 
convergence of 1. Here, the spatially adaptive sparse grid is slightly better than the regular one. 
In $d=3$ dimensions, both sparse grid quadratures are still better than the quasi-Monte Carlo 
quadrature but in the end, the adaptively added points yield worse results. In $d=4$ dimensions, the 
regular sparse grid completely fails, and the spatially adaptive sparse grid is only as good as the 
quasi-Monte Carlo quadrature. In all dimensions, we get the best results with the simplex stochastic 
collocation. In $d=2$ dimensions the maximal obtained order of convergence 2 is twice as good as the 
one for sparse grids and quasi-Monte Carlo quadrature. Additionally, the pre-asymptotic is also 
better. In $d=3$ and $d=4$ dimensions, the convergence rate of the simplex stochastic collocation is 
the same, but, again, the better pre-asymptotic makes a difference. Concluding, we can say that the 
explicit kink approximation is useful and worthwhile, even though the theoretical convergence rates 
are not obtained due to the jumps in the second derivative. All methods requiring a certain 
smoothness suffer from these jumps. 

\section{Conclusion}
We introduced an approach of simplex stochastic collocation for a piecewise approximation of a 
function with polynomials of degree $p$, where the function is not continuously differentiable and 
has kinks. By using a-posteriori information, which incorporates application knowledge, we could 
explicitly approximate the kink, which yields significantly better results. We proved that this 
modification results in algebraic orders of convergence of $(p+1)/d$ and verified the rates with 
test functions in $d=2,3,4$ dimensions. Moreover, we introduced two new error estimators for an 
adaptive refinement. We showed that in contrast to the original error estimators, our ones were 
reliable and solution-based without incorporating unnecessary simulation runs. For multiple 
refinements, as proposed in \cite{witteveen2012_1, witteveen2012_2, witteveen2013}, we analyzed the 
error distribution over the simplices and showed that this approach is reasonable and does not
affect the convergence rates.

We applied our improved version of simplex stochastic collocation to a real gas network. 
Due to the empirical behavior of the employed gas network simulator, which resulted in jumps in the 
second partial derivative of the quantity of interest due to numerical reasons, we could not reach 
the desired convergence rates for the quantity of interest.
Nevertheless, we saw that even in $d=4$ dimensions only 100 sampling points were necessary to 
approximate an expected value as accurate as the model error of $10^{-4}$. A comparison with other 
common methods, such as sparse grid and (quasi-) Monte Carlo quadrature, showed that our method 
benefits from the explicit kink approximation and, hence, yields significantly better results.

So far, we have used the simplex stochastic collocation only for random variables that were 
uniformly distributed. Therefore, the next canonical step will be to extend the method of simplex 
stochastic collocation for random variables following other distributions with bounded support. 
Instead of weighting an error estimator with the area of a simplex, the error estimator could be 
weighted with the probability of a simplex. This idea was already presented for the original version 
of stochastic simplex collocation \cite{witteveen2012_1, witteveen2012_2, witteveen2013} and should 
not cause any problems. The more interesting question is whether simplex stochastic collocation can 
be used for random variables whose density function has unlimited support and how bounding the 
support influences the method. 

Furthermore, we have seen that the method of Voronoi piecewise surrogate models~\cite{rushdi2017} 
provided better convergence results than simplex stochastic collocation for functions with many 
local minima and maxima. This could be due to the fact that, in Voronoi piecewise surrogate models, 
the approximation is based on solving a regression problem over the $2P$-nearest neighbors of each 
cell instead of solving an interpolation problem over the $P$-nearest neighbors. Therefore, it 
should be investigated how the use of regression affects simplex stochastic collocation and whether 
it improves its convergence.

\section*{Acknowledgement}

This work is funded by the German Federal Ministry for Economic Affairs and Energy (BMWi) within the project MathEnergy.

\appendix
\section{Appendix}
\subsection{Generation of Random Points in Simplices}\label{app:a}
According to \cite{devroye1986}, an efficient way to sample uniform distributed random points in the unit 
simplex
\begin{align}
	S_d = \left\{ (s_1, \ldots, s_d): s_i \geq 0, \sum_{i=1}^d s_i \leq 1 \right\}
\end{align}
is the following. Let $u_1, u_2, \ldots, u_{d+1}$ be independent and identically uniform in $[0,1]$ 
distributed random numbers. Then the random variables $e_1=-\log(u_1), e_2 =-\log(u_2),\; \ldots,\; 
\allowbreak e_{d+1}=-\log(u_d)$ are independent and identically exponentially distributed with 
parameter $\lambda=1$. Let $s=\sum_{i=1}^{d-1}s_i$, then the vector 
\begin{align}
	\mathbf{x} = (x_1, x_2, \ldots, x_d) = (e_1/s, e_2/s, \ldots, e_d/s)
\end{align}
is uniform distributed in simplex $S_d$. This method has the advantage that no sample points must be 
rejected nor a sorting of numbers is required.

\bibliographystyle{abbrv}
\bibliography{Literature.bib}

\end{document}